\documentclass[12pt]{amsart}

\usepackage{color}
\definecolor{darkgreen}{rgb}{0, 0.40, 0}

\usepackage[
ps2pdf,                %%% hyper-references for ps2pdf
colorlinks=true,
linkcolor=blue,
citecolor=darkgreen,
urlcolor=magenta
]{hyperref}

\usepackage{epsfig}
\usepackage{latexsym}
\usepackage{psfrag}
\usepackage{amssymb}
\usepackage{amsmath}
\usepackage{amsthm} 
\newcommand{\ie}{{\it i.e.}}

\newcommand{\CC}{{\mathbb{C}}}
\newcommand{\DD}{{\mathbb{D}}}

\renewcommand{\setminus}{{\smallsetminus}}

\newcommand{\lcm}{{\operatorname{lcm}}}

\newcommand{\bigmid}{{~\mid~}}
\newcommand{\st}{{\bigmid}}
\newcommand{\from}{{\colon}} % As in ``f maps _from X _to_ Y''.

\newcommand{\isotopic}{{\medspace \approx \medspace}} % Homeomorphism
\newcommand{\homeo}{{\medspace \cong \medspace}} % Homeomorphism
 
\newcommand{\disjoint}{{\medspace \amalg \medspace}} % Disjointness
\newcommand{\cross}{{\times}}

\newcommand{\lift}[1]{{\widetilde{#1}}}
\newcommand{\closure}[1]{{\overline{#1}}}

\newcommand{\euler}{{\operatorname{\chi}}} % Euler characteristic
\newcommand{\neigh}{{\operatorname{\eta}}} % Open regular neighborhood
\newcommand{\bdy}{{\partial}} % Boundary

\newcommand{\genus}{{\operatorname{genus}}}
\newcommand{\interior}{{\operatorname{interior}}}
\newcommand{\pt}{{\rm pt}} % Point

\theoremstyle{plain}
\newtheorem{theorem}{Theorem}[section]

\newtheorem{lemma}[theorem]{Lemma}
\newtheorem{conjecture}[theorem]{Conjecture}

\theoremstyle{definition}
\newtheorem*{define}{Definition}
\newtheorem{claim}[theorem]{Claim}

\newtheorem{remark}[theorem]{Remark}

\newtheorem*{question}{Question}

\newsavebox{\savepar}

\newcommand{\HcapBdyD}{{c_1}} % |H \cap \bdy D|
\newcommand{\alphacapbeta}{{c_2}} % the maximum of \beta' \cap \alpha'
\newcommand{\lcmnumber}{{c_3}} % the lcm of small positive shifts
\newcommand{\enormous}{{c_4}} % max of the above

\newcommand{\identity}{{\rm Id}}
\newcommand{\last}{{m}}

\begin{document}

\title{Heegaard splittings of the form $H + nK$}
\thanks{This research was supported by Grant No.~2002039 from the
  US-Israel Binational Science Foundation (BSF), Jerusalem, Israel.}
\author{Yoav Moriah} 
\author{Saul Schleimer}
\thanks{The second author was partially supported by NSF Award
  No.~0102069.}
\author{Eric Sedgwick}

\address{\hskip-\parindent
        Yoav Moriah \\
        Department of Mathematics \\
        Technion - Israel Institute of Technology \\
        Haifa 32000, Israel}
\email{ymoriah@tx.technion.ac.il}

\address{\hskip-\parindent
        Saul Schleimer \\
        Department of Mathematics, UIC \\
        851 South Morgan Street \\
        Chicago, Illinois 60607}
\email{saul@math.uic.edu}

\address{\hskip-\parindent
        Eric Sedgwick \\
        DePaul CTI \\
        243 S. Wabash Ave. \\
        Chicago, IL 60604}
\email{esedgwick@cs.depaul.edu}

\date{\today}

\begin{abstract}
Suppose that a three-manifold $M$ contains infinitely many distinct
strongly irreducible Heegaard splittings $H + nK$, obtained by Haken
summing the surface $H$ with $n$ copies of the surface $K$.  We show
that $K$ is incompressible.  All known examples, of manifolds
containing infinitely many irreducible Heegaard splittings, are of
this form.  We also give new examples of such manifolds.
\end{abstract}

\maketitle

\section{Introduction}
\label{Sec:Introduction}

F.~Waldhausen, in his 1978 paper~\cite{Waldhausen78}, asked if every
closed orientable three-manifold contains only finitely many
unstabilized Heegaard splittings.  A.~Casson and C.~Gordon
(see~\cite{CassonGordon85} or~\cite{MoriahSchultens98}), using a
result of R.~Parris~\cite{Parris78}, obtain a definitive ``no''
answer; they obtain examples of closed hyperbolic three-manifolds each
of which contains {\em strongly irreducible} splittings of arbitrarily
large genus.  These examples have been studied and generalized by
T.~Kobayashi~\cite{KobayashiUnpub}, \cite{Kobayashi92}, M.~Lustig and
Y.~Moriah~\cite{LustigMoriah00}, E.~Sedgwick~\cite{Sedgwick97}, and
K.~Hartshorn~\cite{Hartshorn99}.

The goal of this paper is three-fold.  We first show, in
Section~\ref{Sec:OldExamples}, that all of the examples studied so far
are of the form $H + nK$: There is a pair of surfaces $H$ and $K$ in
the manifold so that the strongly irreducible splittings are obtained
via a cut-and-paste construction, {\em Haken sum}, of $H$ with $n$
copies of $K$.  See Section~\ref{Sec:Define} for a precise definition
of Haken sum.

Next, and of more interest, we show when such a sequence exists the
surface $K$ must be incompressible (in Sections~\ref{Sec:HPlusNK}
through~\ref{Sec:HPlusNT}).  We claim:

\begin{theorem}
\label{Thm:HPlusNKGeneral}
Suppose $M$ is a closed, orientable three-manifold and $H$ and $K$ are
closed orientable transverse surfaces in $M$.  Suppose that a Haken
sum $H + K$ is given so that, for arbitrarily large values of $n$, the
surfaces $H + nK$ are pairwise non-isotopic strongly irreducible
Heegaard splittings.  Then the surface $K$ is incompressible.
\end{theorem}

%%% Could make the remark (after the pieces) that this can be made
%%% effective... I.e. we have a bound on how big D (the compressing
%%% disk for K) is -- take D to be fundamental in the triangulation
%%% you produce after being given H and K...

Theorem~\ref{Thm:HPlusNKGeneral} shows that all of the counterexamples
to Waldhausen's question found thus far are Haken manifolds.  This was
already known but required somewhat subtle techniques (see Lemmas~3.2
and~3.3 and Theorem~4.9 of Y.-Q.~Wu's paper~\cite{Wu96}).

Theorem~\ref{Thm:HPlusNKGeneral} was originally conjectured by
Sedgwick along with the much stronger:

\begin{conjecture}
\label{Conj:Sedqwick}
Let $M$ be a closed, orientable $3$-manifold which contains infinitely
many irreducible Heegaard splittings that are pairwise non-isotopic.
Then $M$ is Haken.
\end{conjecture}

We also produce new counterexamples, which are quite different from
those previously studied.  These examples are discussed in
Sections~\ref{Sec:NewExample} through~\ref{Sec:NewExampleIsCorrect}.

The paper concludes in Section~\ref{Sec:Questions} by listing several
conjectures.

\vspace{1mm}
\noindent 
{\bf Acknowledgments:} We thank Tsuyoshi Kobayashi for several
enlightening conversations which led directly to the examples in
Section~\ref{Sec:NewExample}.  We thank David Bachman for bringing the
paper~\cite{KobayashiRieck04} to our attention.  We also like to thank
DePaul University, UIC, and the Technion for their hospitality.

\section{Preliminaries}
\label{Sec:Define}

Fix $M$, a closed, orientable three-manifold.  If $X$ is a submanifold
of $M$ we denote a open regular neighborhood of $X$ by $\neigh(X)$.

A surface $K$ is {\em incompressible} in $M$ if $K$ is embedded,
orientable, closed, not a two-sphere, and a simple closed curve
$\gamma \subset K$ bounds an embedded disk in $M$ if and only if
$\gamma$ bounds a disk in $K$.  The three-manifold $M$ is {\em
  irreducible} if every embedded two-sphere bounds a three-ball in
$M$.  If $M$ is irreducible and contains an incompressible surface
then $M$ is a {\em Haken} manifold.

A surface $H$ is a {\em Heegaard splitting} for $M$ if $H$ is
embedded, connected, and separates $M$ into a pair of handlebodies,
say $V$ and $W$.  A disk $D$ properly embedded in a handlebody $V$ is
{\em essential} if $\bdy D \subset \bdy V$ is not null-homotopic in
$\bdy V$.

\begin{define}
\label{Def:HeegaardSplittingReducible}
A Heegaard splitting $H \subset M$ is {\em reducible} if there is a
pair of essential disks $D \subset V$ and $E \subset W$ with $\bdy D =
\bdy E$.  If $H$ is not reducible it is {\em irreducible}.
\end{define}

\begin{define}
\label{Def:HeegaardSplittingWeaklyReducible}
A Heegaard splitting $H \subset M$ is {\em weakly reducible} if there
is a pair of essential disks $D \subset V$ and $E \subset W$ with
$\bdy D \cap \bdy E = \emptyset$. (See~\cite{CassonGordon87}.)  If $H$
is not weakly reducible it is {\em strongly irreducible}.
\end{define}

One reason to study strongly irreducible Heegaard splittings is that
these surfaces have many of the properties of incompressible surfaces.
An important example of this is:

\begin{lemma}[Scharlemann's No Nesting Lemma~\cite{Scharlemann98}]
\label{Lem:NoNesting}
Suppose that $H \subset M$ is a strongly irreducible Heegaard
splitting and that the simple closed curve $\gamma \subset H$ bounds a
disk $D$ embedded in $M$ and transverse to $H$.  Then $\gamma$ bounds
a disk in either $V$ or $W$. \qed
\end{lemma}

%%% The condition that $D$ be transverse cannot be omitted -- think
%%% about the $(1,1)$ curve on the genus one splitting of $S^3$.
%%% Essentially $D$ and $H$ have to give the same framing to $\gamma$
%%% -- I see the examples but I am not clear on the theory behind this
%%% behavior... 

\begin{figure}[ht]
%%% For figure placement see page 197
\psfrag{sum}{sum}
\psfrag{H}{$F$}
\psfrag{K}{$G$}
\psfrag{H+K}{$F+G$}
\psfrag{H+nK}{$F+nG$}
\psfrag{ncopies}{$n$ copies}
$$\begin{array}{c}
\epsfig{file=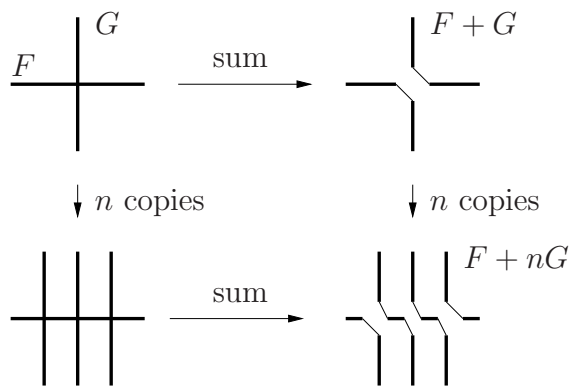, height = 5 cm}
\end{array}$$
\caption{For every intersection of $F$ and $G$ we have $n$
  intersections of $F$ and $nG$.  The light lines are the annuli
  $A_+(\gamma_i)$.} 
\label{Fig:HakenSumExplained}
\end{figure}

We now turn from Heegaard splittings to the concept of the {\em Haken
  sum} of a pair of surfaces.  See Figure~\ref{Fig:HakenSumExplained}
for an illustration.

Suppose $F, G \subset M$ are a pair of closed, orientable, embedded,
transverse surfaces.  Assume that $\Gamma = F \cap G$ is nonempty.
Note that, for every $\gamma \in \Gamma$, the open regular
neighborhood $T(\gamma) = \neigh(\gamma)$ is an open solid torus in
$M$.  Note that $\bdy \closure{T(\gamma)} \setminus (F \cup G)$ is a
union of four open annuli $A_1(\gamma) \cup A_2(\gamma) \cup
A_3(\gamma) \cup A_4(\gamma)$, ordered cyclically.  We collect these
into the two opposite pairs; $A_+(\gamma) = A_1 \cup A_3$ and
$A_-(\gamma) = A_2 \cup A_4$.  For every $\gamma \in \Gamma$ we now
chose an $\epsilon(\gamma) \in \{+, -\}$ and form the {\em Haken sum}:

$$F + G = \left( \left( F \cup G \right) \setminus \left(
\bigcup_\gamma T(\gamma) \right) \right) \cup \left( \bigcup_\gamma
A_{\epsilon(\gamma)}(\gamma) \right)$$

Note that the Haken sum depends heavily on our choices of
$\epsilon(\gamma)$.  As a bit of notation we call the core curves of
the annuli $A_\epsilon$ the {\em seams} of the Haken sum.  Also there
is an obvious generalization of Haken sum to properly embedded
surfaces.

\begin{remark}
\label{Rem:HakenSumNotCanonical}
If $F$ and $G$ are compatible normal surfaces, carried by a single
branched surface, or transversely oriented there is a natural choice
for the function $\epsilon(\gamma)$.  
\end{remark}

We now define the Haken sum $F + nG$: Take $n$ parallel copies of $G$
in $\neigh(G)$ and number these $\{G_i\}_1^n$.  For every curve
$\gamma \in \Gamma$ we now have $n$ curves $\{\gamma_i \subset F \cap
G_i\}_{i=1}^n$.  A Haken sum $F + G$ is determined by labellings
$A_\pm(\gamma)$ and choices $\epsilon(\gamma) \in \{+, -\}$.  Using
the parallelism of the $G_i$ we take identical labellings for
$A_\pm(\gamma_i)$ and make identical choices for $\epsilon(\gamma_i)$.
See Figure~\ref{Fig:HakenSumExplained} for a cross-sectional view at
$\gamma$.

The surface $F + nG$ is now the usual Haken sum of $F$ and $nG$ with
these induced choices, $A_\pm(\gamma_i)$ and $\epsilon(\gamma_i)$.

\section{Existing examples}
\label{Sec:OldExamples}

This section shows that the Casson-Gordon examples are of the
form $H + nK$.  At the end of the section we briefly discuss the
examples of Kobayashi~\cite{Kobayashi92}, and
Lustig and Moriah~\cite{LustigMoriah00}.

\begin{figure}[ht]
\psfrag{k}{$k$}
\psfrag{S}{$S$}
$$\begin{array}{c}
\epsfig{file=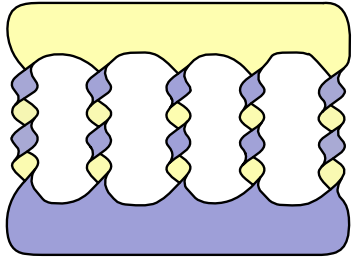, height = 3.5 cm}
\end{array}$$
\caption{The $k(5,5,5,5,5)$-pretzel knot.}
\label{Fig:Pretzel}
\end{figure}

Let $k = k(n_1, \ldots, n_m) \subset S^3$ be a {\em pretzel
  knot}~\cite{Kawauchi85} with {\em twist boxes} of order $n_i$.  Here
we choose $m$ and the $n_i$ to be odd, positive, and greater than $4$.
See Figure~\ref{Fig:Pretzel} for an example.

\begin{figure}[ht]
\psfrag{k}{$k_1$}
\psfrag{S}{$S$}
$$\begin{array}{c}
\epsfig{file=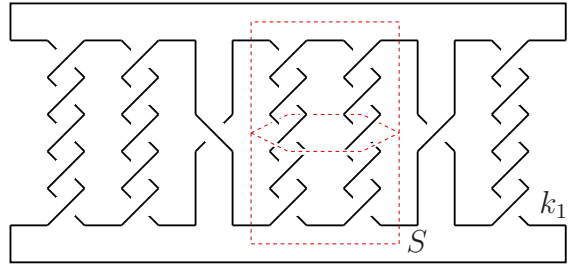, height = 3.5 cm}
\end{array}$$
\caption{After twisting the $k(5,5,5,5,5)$-pretzel knot we obtain the
 $k(5,5,-1,5,5,1,5)$-pretzel knot.}
\label{Fig:TwistedPretzel}
\end{figure}

A pretzel knot has an associated Seifert surface, $F$.  This is the
compact checkerboard surface for the standard diagram. Again, see
Figure~\ref{Fig:Pretzel}.  Let $B$ be the three-ball containing the
pair of consecutive twist boxes of order $n_i$ and $n_{i+1}$.  Let $S
= \bdy B$.  Note that $|k \cap S| = 4$; see
Figure~\ref{Fig:TwistedPretzel}.  There is a well-known twisting
procedure which twists $k = k(n_1, \ldots, n_m)$ along $S$ giving 
$$k_1
= k(n_1, \ldots, n_{i-1}, -1, n_i, n_{i+1}, 1, n_{i+2}, \ldots, n_m).$$
Again, see Figure~\ref{Fig:TwistedPretzel}.

So, given the pretzel knot $k$ and the sphere $S$ we can produce the
sequence $\{k_n\}$ of {\em $n$-times twisted} pretzels:

$$k_n = k(n_1, \ldots, n_{i-1} , \overbrace{-1, \ldots, -1}^n, n_{i},
n_{i+1}, \overbrace{1, \ldots, 1}^n, n_{i+2}, \ldots, n_m).$$

\noindent Denote the associated Seifert surface for $k_n$ by
$F_n$.  Note that $k_n$ is isotopic to $k = k_0$ and that $F_0 = F$.

In his thesis, Parris proves:

\begin{theorem}[Parris~\cite{Parris78}]
\label{Thm:Parris}
The surfaces $F_n$ are free incompressible Seifert surfaces for
$k$. \qed
\end{theorem}

Let $X = S^3 \setminus \neigh(k_n)$.  Let $\widehat{V_n}$ be a closed
regular neighborhood of $F_n \cup \neigh(k_n)$.  So $k_n \subset
\widehat{V_n}$.  Let $W_n = \closure{S^3 \setminus \widehat{V_n}}$.
Now, as $k_n$ is isotopic into $H_n = \bdy \widehat{V_n}$, doing $1/l$
Dehn surgery along $k$ makes $\widehat{V_n}$ into a handlebody, which
we denote by $V_n$.  Here $l$ is any positive integer greater than
$4$.  Let $M = X(1/l)$ be the $1/l$ Dehn surgery of $S^3$ along $k$.
Let $H_n = \bdy V_n = \bdy W_n \subset M$.  Note that the genus of
$H_n$ is $2n + 4$.  We have:

\begin{theorem}[Casson and
  Gordon~\cite{CassonGordon85},~\cite{MoriahSchultens98}]
\label{Thm:CassonGordon}
The Heegaard splittings $H_n \subset M$ are strongly irreducible. \qed
\end{theorem}

Now let $G$ be the surface $\bdy(B \setminus \neigh(k)) = (S \setminus
\neigh(k)) \cup (\bdy\closure{\neigh(k)} \cap B)$.  We now state the
main theorem of this section:

\begin{theorem}
\label{Thm:CGPAreHPlusNK}
The Heegaard surfaces $H_n$ are isotopic to a Haken sum $H_0 + 2nG$.
\end{theorem}

We require several lemmas for the proof of
Theorem~\ref{Thm:CGPAreHPlusNK}.

\begin{lemma}
\label{Lem:FPlusS}
The surface $F_n$ is isotopic to $F_0 + nS$.
\end{lemma}

\begin{figure}[ht]
$$\begin{array}{ccc}
\epsfig{file=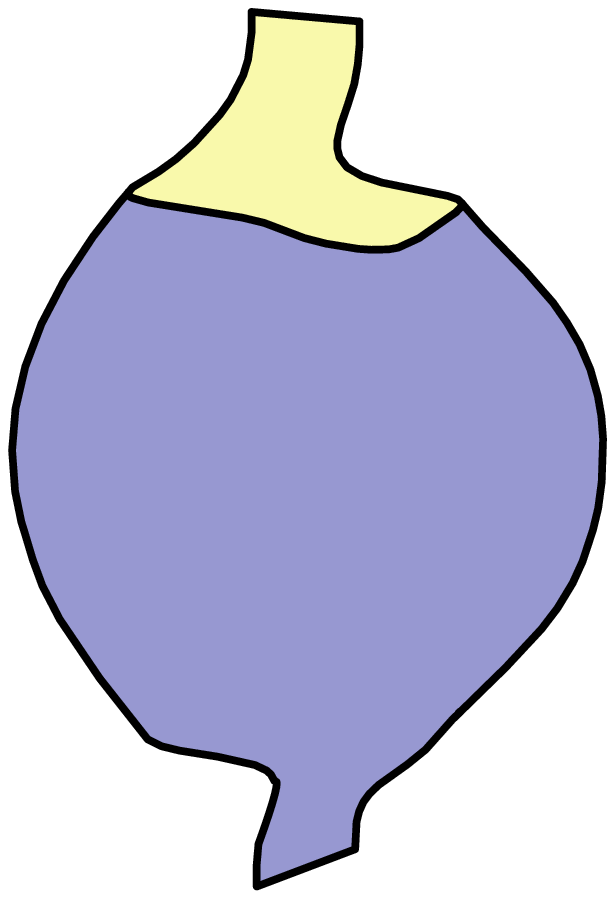, height = 3.5 cm} &
\epsfig{file=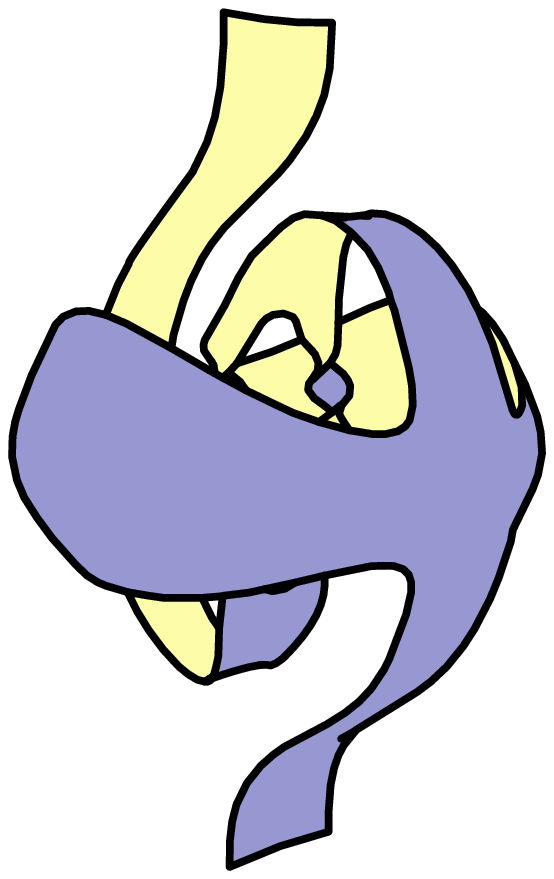, height = 3.5 cm} &
\epsfig{file=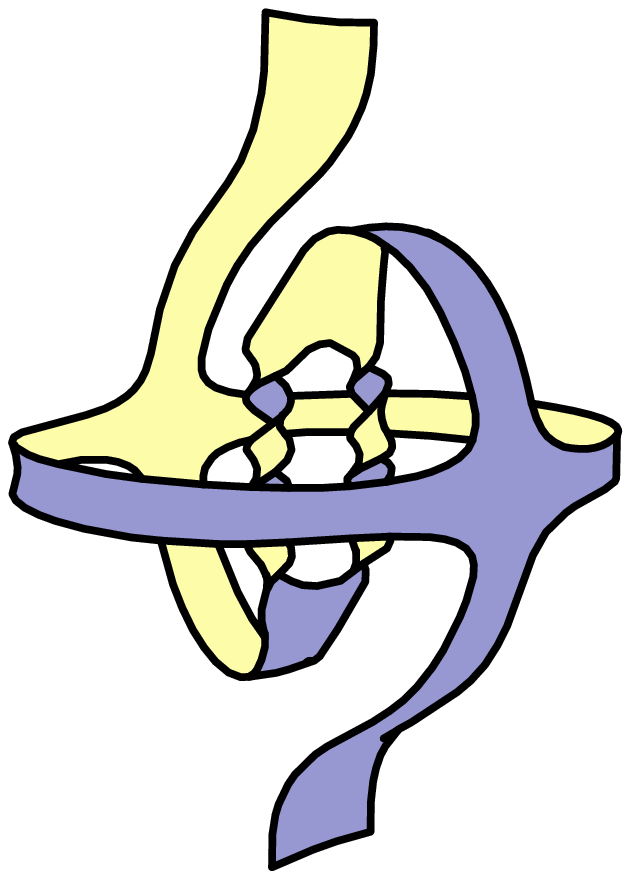, height = 3.5 cm} \\
\end{array}$$
\caption{Isotoping $F + S$, moving $k$ near the equator of $S$.}
\label{Fig:IsotopyOfF+S}
\end{figure}

\begin{figure}[ht]
$$\begin{array}{c}
\epsfig{file=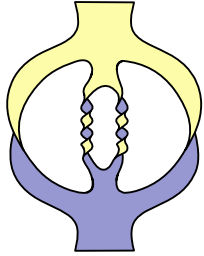, height = 3.5 cm}
\end{array}$$
\caption{Flatten the resulting figure into the plane of the diagram.}
\label{Fig:Flatten}
\end{figure}

\begin{proof}
Let $\alpha$ and $\beta$ be the arcs of intersection between $S$ and
$F = F_0$.  Let $B_\alpha$ be a closed regular neighborhood of
$\alpha$.  Let $S_\alpha$ be the boundary of $B_\alpha$.  See the left
side of Figure~\ref{Fig:B_Alpha} for a picture of $S \cup F$ inside
$B_\alpha$.

\begin{figure}[ht]
\psfrag{K}{$k$}
\psfrag{S}{$S$}
\psfrag{F}{$F$}
$$\begin{array}{ccc}
\epsfig{file=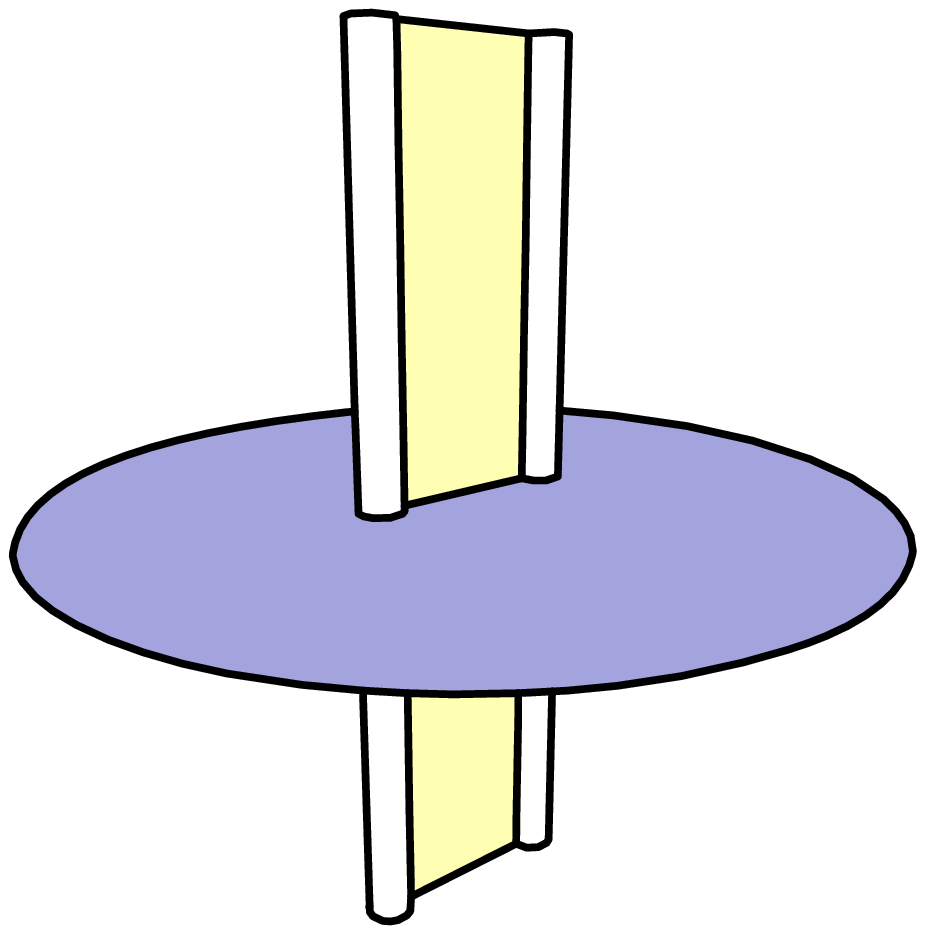, height = 3.5 cm} &
\epsfig{file=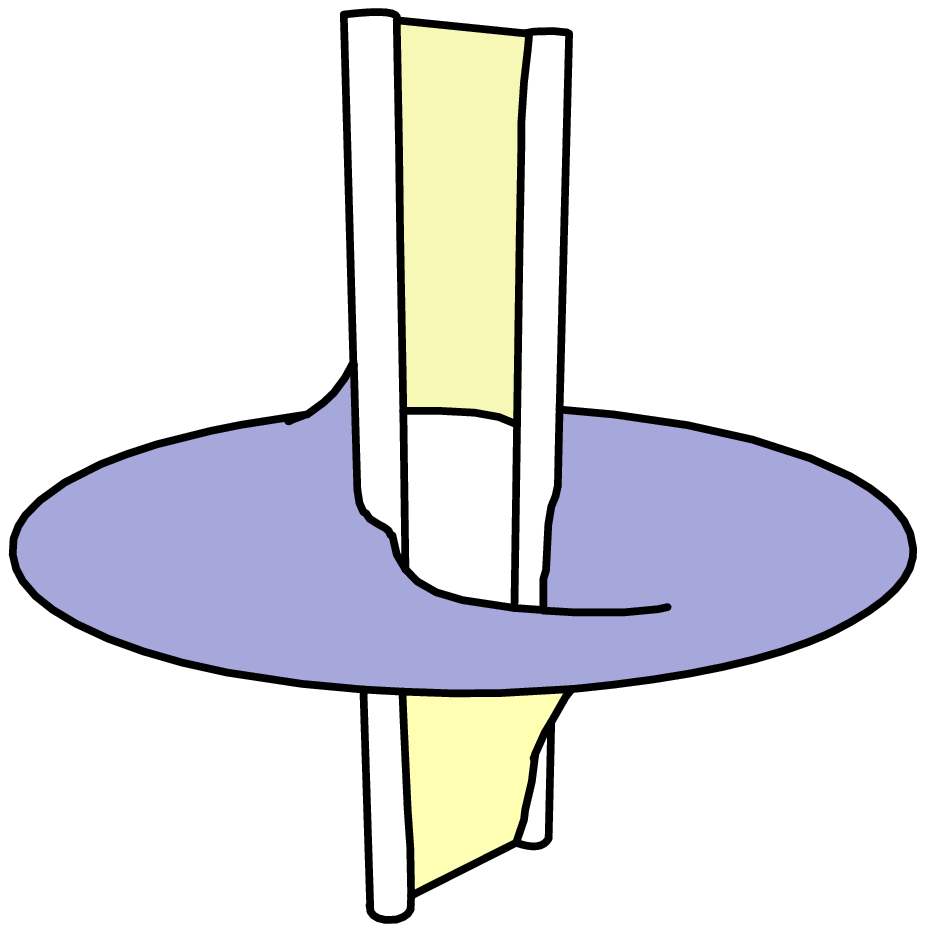, height = 3.5 cm} &
\epsfig{file=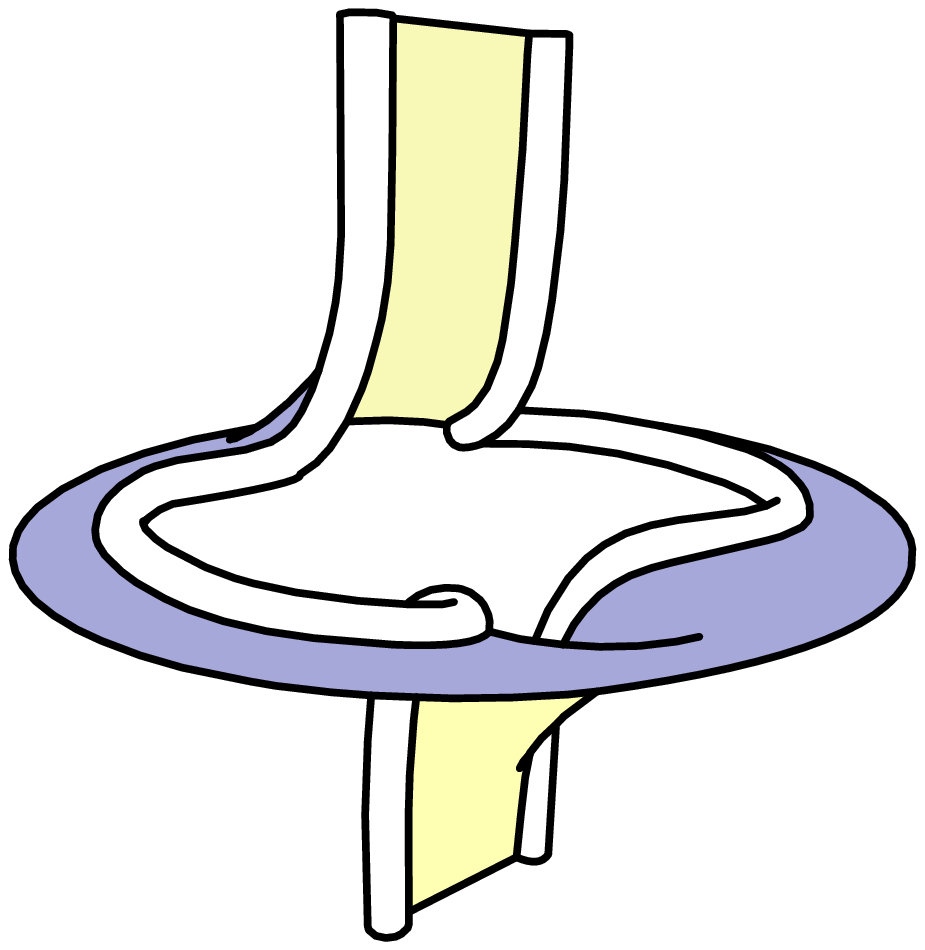, height = 3.5 cm}
\end{array}$$
\caption{The knot $k$ has been thickened a bit.  On the left $F$ is
  vertical while $S$ is horizontal.  The middle is their Haken sum.  The
  right shows the isotopy of $\alpha' \cup k' \cup \alpha'' \cup k''$
  to be horizontal.}
\label{Fig:B_Alpha}
\end{figure}

We choose the Haken sum which glues the top sheet of $(F \cap
B_\alpha) \setminus \alpha$ to the back sheet of $(S \cap B_\alpha)
\setminus \alpha$.  Glue the bottom sheet of $(F \cap B_\alpha)
\setminus \alpha$ to the front sheet of $(S \cap B_\alpha) \setminus
\alpha$.  See the right hand side of Figure~\ref{Fig:B_Alpha} for a
picture of the Haken sum.

Let $\alpha'$ and $\alpha''$ be the seams along which the sheets of
$F$ and $S$ are glued.  Let $k'$ and $k''$ be the arcs of $k \setminus
(\bdy \alpha' \cup \bdy \alpha'')$ inside of $B_\alpha$.  Do a small
isotopy of the loop $\gamma = \alpha' \cup k' \cup \alpha'' \cup k''$
as shown in Figure~\ref{Fig:B_Alpha}.  After this isotopy the image of
$\gamma$ lies in a regular neighborhood of the curve $S_\alpha \cap
S$.

We perform the same sequence of steps near $\beta$.  Recall that
$S_\alpha \cap S$ and $S_\beta \cap S$ cobound an annulus, $A \subset
S$. Isotope the surface $F + S$ to move $k$ close to the core curve of
$A$ -- this isotopy is illustrated in a sequence of steps in
Figure~\ref{Fig:IsotopyOfF+S}.

Now flatten out the right hand side of Figure~\ref{Fig:IsotopyOfF+S}
by rotating the two twist boxes inside of $S$ by $180^\circ$.  Also
flatten the annulus into the plane containing the standard diagram of
$k$.  See Figure~\ref{Fig:Flatten}.

Note that the result is the Seifert surface associated to the pretzel
knot $k_1 = k(5, 5, -1, 5, 5, 1, 5)$.  Thus, by induction, the proof of
Lemma~\ref{Lem:FPlusS} is complete.
\end{proof}

Recall that $k$ is the given pretzel knot, $F = F_0$ is the associated
Seifert surface, and $S$ is the two-sphere bounding the three-ball
$B$, as above.

\begin{lemma}
\label{Lem:FPlusG}
The surface $F_n$ is isotopic to $F_0 + nG$.
\end{lemma}

\begin{figure}[ht]
$$\begin{array}{cc}
\epsfig{file=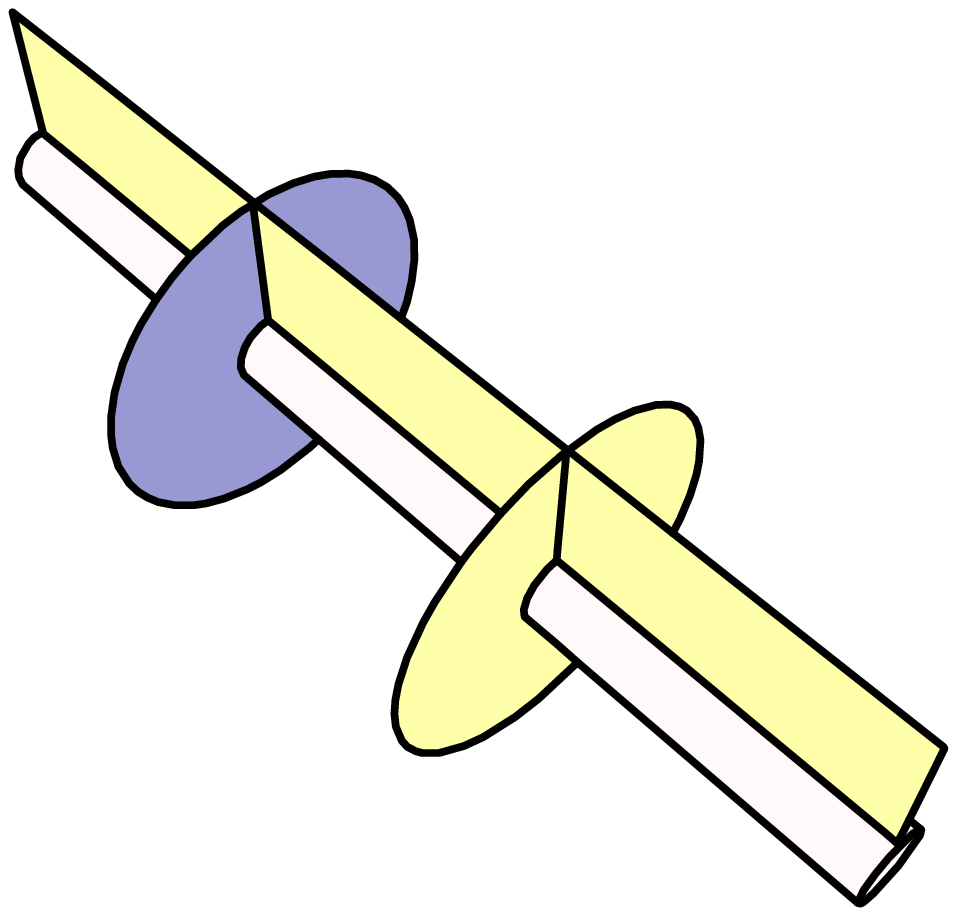, height = 3.5 cm} &
\epsfig{file=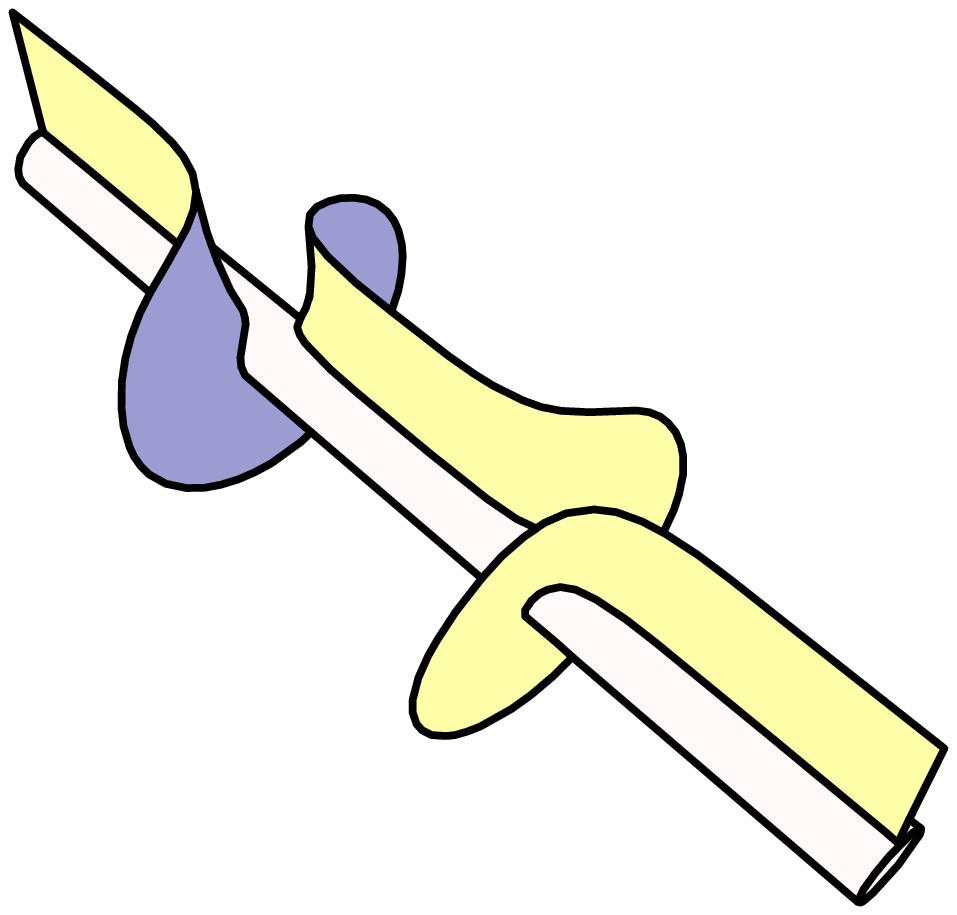, height = 3.5 cm} 
\end{array}$$
\caption{Forming the Haken sum of $F$ (longitudinal) and $S$ (meridional).}
\label{Fig:FAndSNearK}
\end{figure}

\begin{figure}[ht]
$$\begin{array}{cc}
\epsfig{file=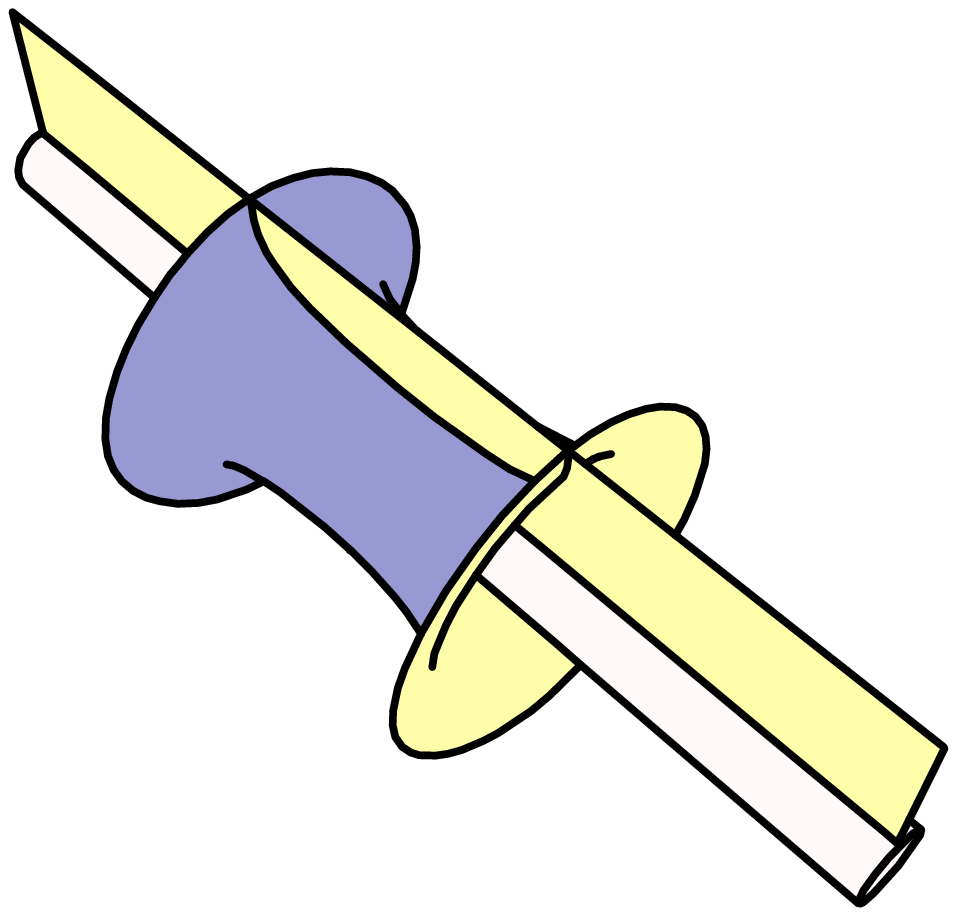, height = 3.5 cm} &
\epsfig{file=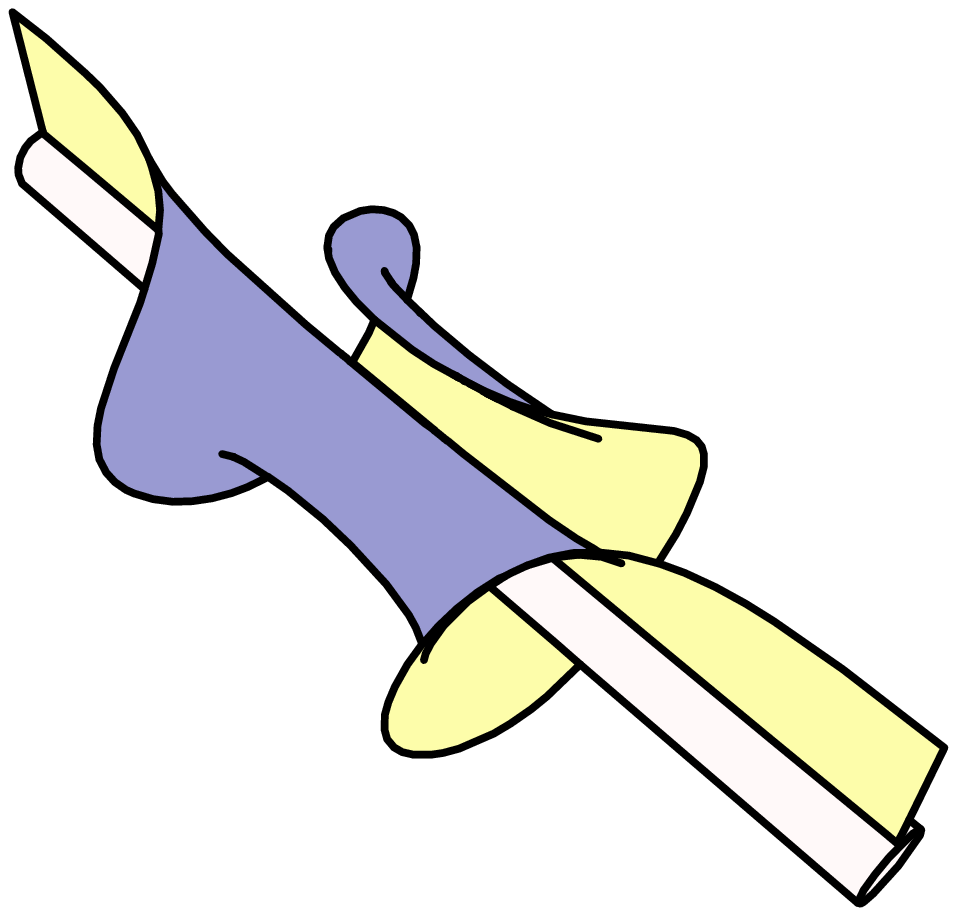, height = 3.5 cm}
\end{array}$$
\caption{Forming the Haken sum of $F$ and $G$.}
\label{Fig:FAndGNearK}
\end{figure}

\begin{proof}
Consider a single component of $\neigh(k) \cap \neigh(B)$.  This
component $B'$ is a ball.  Let $k' = k \cap B'$.  The disk $F' = F
\cap B'$ is a boundary compression for $k'$ in $B'$.  The two disks
$S' \cup S'' = S \cap B'$ each intersect $k'$ in a single point.  See
the left hand side of Figure~\ref{Fig:FAndSNearK} for a picture.  (The
knot $k$ has been thickened a bit.)

The arcs $S' \cap F'$ and $S'' \cap F'$ are both part of $\alpha
\subset S \cap F$.  Thus Haken summing along $S' \cap F'$ agrees with
Haken summing along $S'' \cap F'$.  See the right hand side of
Figure~\ref{Fig:FAndSNearK}.

Turn now to $F + G$.  Recall that $G = \bdy(B \setminus \neigh(k))$.
Note that $G \cap \closure{\neigh(k)}$ is a pair of annuli.  Isotope
these annuli, rel boundary, slightly into $\neigh(k)$ so that $G
\setminus \neigh(k)$ is identical to $S \setminus \neigh(k)$.  Thus
obtain the picture of $F \cap B'$ and $G \cap B'$ shown on the left in
Figure~\ref{Fig:FAndGNearK}.

Finally take the Haken sum of $F'$ with $G' = G \cap B'$ as forced by
our previous choices.  See the right of Figure~\ref{Fig:FAndGNearK}.
Note that $F' + G'$ is isotopic to $F' + (S' \cup S'')$, rel boundary.
The same holds inside the other component of $\neigh(k) \cap B$.
Finally $F + S$ is identical to $F + G$ outside of $\neigh(k)$.  The
lemma is proved.
\end{proof}

We are now equipped to prove Theorem~\ref{Thm:CGPAreHPlusNK}:

\begin{proof}
Notice now that $H_n$ is isotopic to the boundary of a regular
neighborhood of $F_n$.  As $\bdy F_n = k_n$ the splitting $H_n$ is
obtained by gluing two parallel copies of $F_n$ with an annulus $A_n
\subset \bdy\closure{\neigh(k_n)}$, where the core curve of $A_n$ has
longitudinal slope $\bdy\closure{\neigh(k_n)} \cap F_n$.  Note that
$A_0$ is taken to $A_n$ by the twisting isotopy taking $k = k_0$ to
$k_n$. We thus have the following:

\begin{eqnarray}
\label{Eqn:CGP}
H_n & = & 2F_n \cup A_n \\
    & \isotopic & 2(F_0 + nG) \cup A_0 \\
    & = & (2F_0 \cup A_0) + 2nG \\
    & = & H_0 + 2nG.
\end{eqnarray}

The second line follows from Lemma~\ref{Lem:FPlusG}.  The third line
holds because $G$ has no boundary.  This concludes the proof of
Theorem~\ref{Thm:CGPAreHPlusNK}.
\end{proof}

\begin{remark}
\label{Rem:OtherExamples}
The examples of~\cite{Kobayashi92} and~\cite{LustigMoriah00} are very
similar -- they begin with a knot admitting a Conway sphere $S$ and a
natural Seifert surface $F$.  They then isotope the knot by twisting
inside $S$.  Thus their examples of high genus Heegaard splittings may
also be obtained via Haken sum.
\end{remark}

\section{Removing trivial curves}

Here we discuss a method for ``cleaning'' Haken sums.  To be
precise, we have:

\begin{lemma}
\label{Lem:RemoveTrivial}
Suppose $H + nK$ is a sequence of Haken sums.  Let $m$ be the number
of curves of $H \cap K$ which are inessential on $K$.  Then there is
an isotopy of $H' = H + mK$ and a Haken sum $H' + K$ so that
\begin{itemize}
\item
all curves of $H' \cap K$ are essential on $K$ and
\item
for all $n > m$ the surface $H + nK$ is isotopic to $H' + (n-m)K$.
\end{itemize}
\end{lemma}

We call such sequences {\em essential} in $K$.

\begin{proof}[Proof of Lemma~\ref{Lem:RemoveTrivial}]
If $m = 0$ there is nothing to prove.  If not, we claim there is a
surface $\widehat{H}$ such that $\widehat{H}$ is isotopic to $H + K$,
$\widehat{H} \cap K$ has fewer inessential (on $K$) curves than $H
\cap K$ does, and $\widehat{H} + (n - 1)K$ is isotopic to $H + nK$ for
all $n > 0$.  Applying this $m$ times will prove the lemma.

So suppose $\alpha \subset H \cap K$ is inessential on $K$.  Assume that
the disk $D \subset K$ bounded by $\alpha$ is {\em innermost}.  That
is, $D \cap H = \alpha$.

Let $N = \closure{\neigh(K)} \homeo K \cross [0,1]$.  We identify $K$
with $K \cross \{1/2\}$.  Let $D'$ be the component of $(H + K)
\setminus \bdy N$ containing $D$.  Suppose that $\closure{D'}$ has
boundary in $K \cross \{1\}$.  (The case $\closure{D'} \subset K
\cross \{0\}$ is similar.)

Isotope $D'$ up, relative to $(H + K) \cap \bdy N$, to lie in
$\neigh(K \cross \{1\})$, while isotoping all other components of $K
\setminus H$ down into $\neigh(K \cross \{0\})$.  See
Figure~\ref{Fig:RemoveTrivial}.

\begin{figure}[ht]
\psfrag{H+K}{$H + K$}
\psfrag{D}{$D$}
\psfrag{H'}{$\widehat{H}$}
\psfrag{Kx1}{$K \cross \{1\}$}
\psfrag{Kx0}{$K \cross \{0\}$}
$$\begin{array}{c}
\epsfig{file=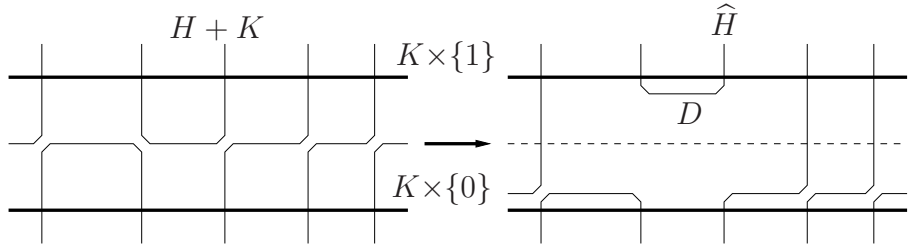, height = 3 cm}
\end{array}$$
\caption{On the left we see $H + K$ intersecting
$\closure{\neigh(K)}$.  On the right $H + K$ has been isotoped to
$\widehat{H}$.}
\label{Fig:RemoveTrivial}
\end{figure}

Let $\widehat{H}$ be this new position of $H + K$ and note that
$\widehat{H} \cap (K \cross \{1/2\})$ has at least one fewer trivial
curve of intersection with $K$.  

We now must prove that $\widehat{H} + (n - 1)K$ is isotopic to $H +
nK$, for all $n > 0$.  Recall that $\alpha$ was the chosen innermost
curve of $H \cap K$, bounding $D \subset K$.  Form $H + nK$ and
isotope all subdisks parallel to $D$ up.  Isotope lowest copy of $K
\setminus D$ down.  This yields $\widehat{H} + (n - 1)K$.  (See
Figure~\ref{Fig:IsotopeTrivial}.)  This completes the claim and thus
the lemma.

\begin{figure}[ht]
\psfrag{H+3K}{$H + 3K$}
\psfrag{H'+2K}{$\widehat{H} + 2K$}
$$\begin{array}{c}
\epsfig{file=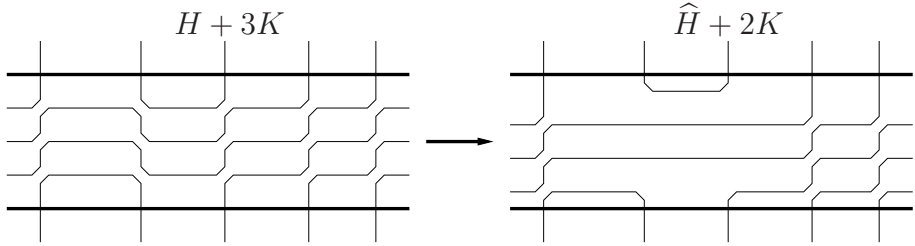, height = 3 cm}
\end{array}$$
\caption{$H + 3K$ is isotopic to $\widehat{H} + 2K$.}
\label{Fig:IsotopeTrivial}
\end{figure}
\end{proof}

\section{Adding surfaces of genus greater than two}
\label{Sec:HPlusNK}

Theorem~\ref{Thm:HPlusNKGeneral} divides into two statements.  The
first addresses the case $\genus(K) > 1$ while the second deals with
the case $K$ a torus.  We begin with:

\begin{theorem}
\label{Thm:HPlusNK}
Suppose $M$ is a closed, orientable three-manifold and $H$ and $K$ are
closed orientable transverse surfaces in $M$, with $\genus(K) \geq 2$.
Suppose that a Haken sum $H + K$ is given so that the surface $H + nK$
is a strongly irreducible Heegaard splitting for arbitrarily large
values of $n$.  Then the surface $K$ is incompressible.
\end{theorem}

We begin by giving a brief sketch of the proof.  Aiming for a
contradiction we assume that $K$ is compressible.  Using
Lemma~\ref{Lem:SeparatingCompression} below we find a compressing disk
$D$ for $K$ with $\bdy D$ separating in $K$.  

For large $n$ the disk $D$ intersects $H + nK$ in a fairly controlled
way -- in particular there is a large family of parallel curves
$\{\gamma_i\}$ in the intersection $(H + nK) \cap D$.  We will show
that many of the $\{\gamma_i\}$ are essential curves on $H + nK$.  By
Scharlemann's ``No Nesting'' Lemma~\ref{Lem:NoNesting} all of these
$\gamma_i$'s bound disks $D_i$ in one of the two handlebodies $V_n$ or
$W_n$.  (Here $\bdy V_n = \bdy W_n$ equals $H + nK$.)  Finally the two
curves $\gamma_i$ and $\gamma_{i+1}$ cobound a subannulus $A_i \subset
D$.  Compressing or boundary compressing $A_i$ will give an essential
disk $E_i$ disjoint from $D_i$.  This demonstrates that $H + nK$
is weakly reducible, a contradiction.

\pagebreak
%%% See page 214 of latex book for discussion of pagebreaks.

\subsection{Finding a separating compressing disk}

We will need a simple lemma:
\begin{lemma}
\label{Lem:SeparatingCompression}
If $G \subset M$ is a compressible surface, which is not a torus, then
there is a compressing disk $D \subset M$ so that $\bdy D$ is a
separating curve on $G$.
\end{lemma}

\begin{proof}
Let $E$ be any compressing disk for $G$.  If $\bdy E$ is a separating
curve then take $D = E$ and we are done.  So suppose instead that
$\bdy E$ is non-separating in $G$.  Choose $\gamma \subset G$ to be
any simple closed curve which meets $\bdy E$ exactly once.  Let $N$ be
a closed regular neighborhood of $\gamma \cup E$, taken in $M$.  Let
$D$ be the closure of the disk component of $\bdy N \setminus G$.
This is the desired disk.
\end{proof}

\subsection{The intersection with the compressing disk}
\label{Sec:Assumptions}

We now begin the proof of Theorem~\ref{Thm:HPlusNK}.  

Recall that $H$ and $K$ are a pair of surfaces so that $H + nK$ is a
strongly irreducible Heegaard splitting for arbitrarily large $n$.
Applying Lemma~\ref{Lem:RemoveTrivial} we may assume that every curve
of intersection between $H$ and $K$ is essential in $K$.

In order to obtain a contradiction assume that $K$ is compressible.
Use Lemma~\ref{Lem:SeparatingCompression} to obtain a compressing disk
$D$ for $K$, transverse to $H$, where $\bdy D$ is separating in $K$.
We may choose $D$ to minimize the size of the intersection $|(H \cap
K) \cap D|$.  Denote the two components of $K \setminus \bdy D$ by
$K'$ and $K''$.

\begin{figure}[ht]
\psfrag{stack}{stack}
\psfrag{H}{$H$}
\psfrag{K1}{$K_1$}
\psfrag{Kn}{$K_4$}
\psfrag{bdyD}{$\bdy D$}
$$\begin{array}{c}
\epsfig{file=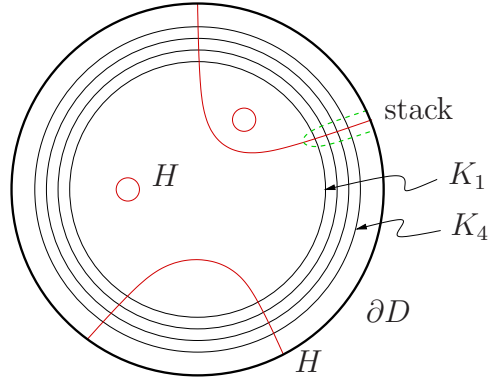, height = 5 cm}
\end{array}$$
\caption{A picture of $D$.  The concentric circles are the curves of
  $K_i \cap D$.  The arcs and small circles make up $H \cap D$.}
\label{Fig:PictureOfD}
\end{figure}

For any $n > 0$ such that $H + nK$ is a strongly irreducible Heegaard
splitting proceed as follows: Label the components of $nK$ by $K_1,
\ldots, K_n$. Isotopy $nK$ so that all of the $K_i$ lie inside of
$\neigh(K)$, are disjoint from $K$, and meet $\interior(D)$ in a
single curve.  Choose subscripts for the $K_i$ consecutively so that
$K_1 \cap D$ is innermost among the curves of intersection $(\cup K_i)
\cap D$.  See Figure~\ref{Fig:PictureOfD} for a picture of how the
$K_i$ and $H$ intersect $D$.

Note $H \cap D$ is a collection of arcs and simple closed curves.
The arcs' intersection with $K_i \cap D$ will give a cross-sectional
view of the Haken sum of $H$ with $nK$.

\begin{figure}[ht]
\psfrag{bdyD}{$\bdy D$}
\psfrag{H+nK}{$H+nK$}
\psfrag{plus}{plus}
\psfrag{minus}{minus}
$$\begin{array}{cc}
\epsfig{file=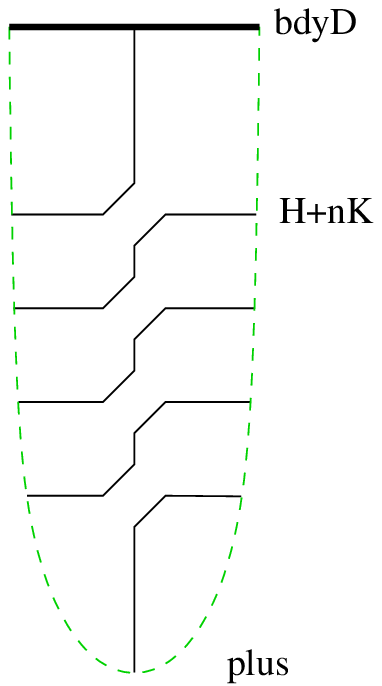, height = 3.5 cm} &
\epsfig{file=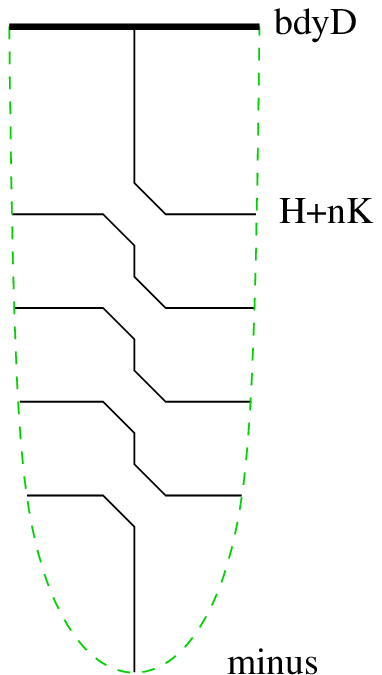, height = 3.5 cm}
\end{array}$$
\caption{In both cases we are looking at $D$ from the direction which
  the transverse orientation points in.}
\label{Fig:Stacked}
\end{figure}

Fix attention on a {\em stack} of intersections, \ie, a collection of $n$
consecutive points of intersection between an arc of $H \cap D$ and
$nK$, all of which are close to a point of $H \cap \bdy D$.  Again,
see Figure~\ref{Fig:PictureOfD}.  Choose a transverse orientation on
$D$.  Assign a parity to the stack as follows: A stack is {\em
  positive} if, after the Haken sum, the segment of $(K_i \cap D)
\setminus \neigh(K_i \cap H)$ on the left is attached to the segment
of $(K_{i+1} \cap D) \setminus \neigh(K_{i+1} \cap H)$ on the right.
Otherwise the stack is {\em negative}.  See Figure~\ref{Fig:Stacked}.

\begin{claim}
\label{Clm:Stacks}
The number of positive stacks equals the number of negative stacks.
\end{claim}

\begin{proof}
Recall $\bdy D$ separates $K$ into two pieces, $K'$ and $K''$.  So
every component of $H \cap K'$ is either a simple closed curve,
disjoint from $\bdy D$, or is a properly embedded arc.  Pick one of
these arcs, say $\alpha \subset H \cap K'$.  Note the endpoints of
$\alpha$ lie in $\bdy D$ and give rise to stacks of opposite parity.
\end{proof}

\begin{figure}[ht]
\psfrag{bdyD}{$\bdy D$}
\psfrag{gamma1}{$\gamma_3$}
\psfrag{H+nK}{$H + nK$}
\psfrag{x}{$x$}
\psfrag{x1}{$x_1$}
\psfrag{x4}{$x_4$}
$$\begin{array}{c}
\epsfig{file=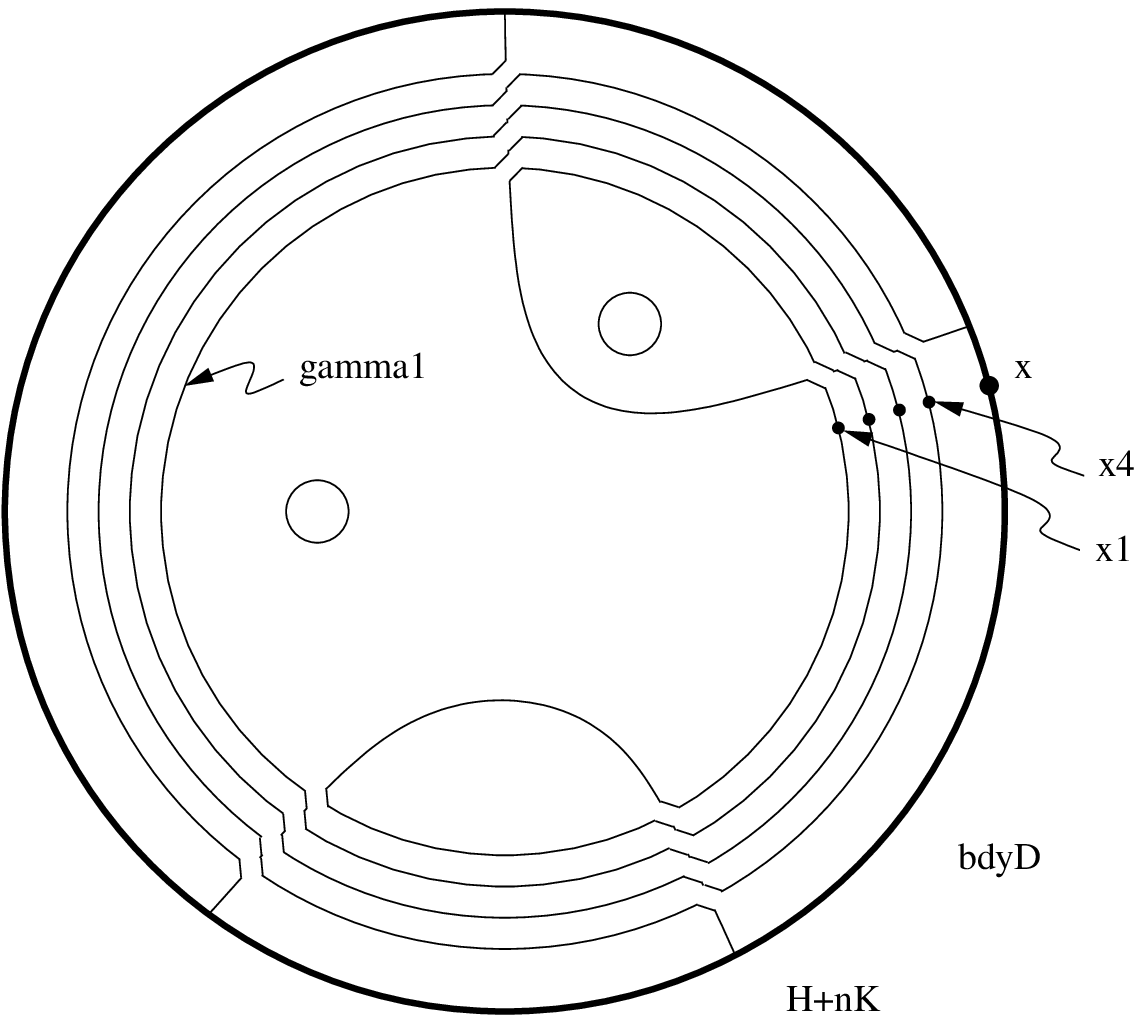, height = 5 cm}
\end{array}$$
\caption{}
\label{Fig:Target}
\end{figure}

Next, analyze how the intersection $(H + nK) \cap D$ lies in $D$:  As in
Figure~\ref{Fig:Target} fix any point $x \in (\bdy D \setminus H)$.
Let $x_i$ be the corresponding point of $K_i \cap D$.

An arc of $(K_i \cap D) \setminus \neigh(H \cap nK)$ is a {\em
  horizontal arc at level $i$}.  In particular the arc containing
$x_i$ is at level $i$.  Orient these arcs in a clockwise fashion.
Note that horizontal arcs are also subarcs of $(H + nK) \cap D$.  When
a horizontal arc at level $i$ enters a positive stack it {\em ascends}
and when it enters a negative stack it {\em descends} a single level.

Consider now an arc of $(H \cap D) \setminus \neigh(H \cap nK)$.
These are the {\em vertical} arcs.  If a vertical arc meets $\bdy D$
call it an {\em external} arc.  If a vertical arc is contained in the
subdisk of $D$ bounded by $K_1 \cap D$ call it an {\em internal} arc.
See Figure~\ref{Fig:PictureOfD}.

Suppose the component of $(H + nK) \cap D$ which contains $x_i$ does
{\em not} contain any internal or external vertical arcs.  Then call
that component $\gamma_i$.  For each value of $i$ where the property
above does not hold, $\gamma_i$ is left undefined.

Set 
\begin{equation}
\label{Eqn:HcapBdyD}
\HcapBdyD = |H \cap \bdy D|.
\end{equation}
Note that $\HcapBdyD$ is even.

\begin{claim}
\label{Clm:ArcsAndCurves}
The collection $(H + nK) \cap D$ consists of
\begin{itemize}
\item
exactly $\HcapBdyD/2$ arcs, 
\item
the curves $\{\gamma_i\}$, and
\item
at most another $|H \cap D|$ simple closed curves.
\end{itemize}
Furthermore, each $\gamma_i$ is a simple closed curve.  Also
$|\{\gamma_i\}| \geq n - \HcapBdyD$.  Finally, $\gamma_i$ and
$\gamma_{i+1}$ cobound an annulus component $A_i$ of $D \setminus (H +
nK)$.
\end{claim}

The claim follows from Figure~\ref{Fig:Target}.  For completeness, a
proof is included.

\begin{proof}[Proof of Claim~\ref{Clm:ArcsAndCurves}.]
The first statement in the claim is trivial: $H \cap \bdy D$ and $(H
+ nK) \cap \bdy D$ are the same set of points.  Next, count the
$\gamma_i$'s:

Choose any $i$ with $\HcapBdyD/2 < i < n - \HcapBdyD/2$ and let
$\alpha$ be the component of $(H + nK) \cap D$ containing $x_i$.
Starting at $x_i$, and moving along $\alpha$ in a clockwise fashion,
we ascend whenever we go through a positive stack and descend through
the negative stacks.  As there are $\HcapBdyD/2$ positive stacks and
the same number of negative stacks $\alpha$ contains no internal or
external vertical arcs.  Also $\alpha$ goes through none of the other
$x_j$'s.  So $\alpha$ is a simple closed curve and is labelled
$\gamma_i$.

It follows there are at least $n - \HcapBdyD$ of the
$\gamma_i$'s in $(H + nK) \cap D$.  These are all parallel in $D$,
yielding the annuli $\{A_i\}$.  Again, see Figure~\ref{Fig:Target}.

To finish the claim note that any simple closed curve of $(H + nK)
\cap D$, which is not a $\gamma_i$, is either a simple closed curve
component of $H \cap D$ or contains an internal vertical arc.  Thus
there are at most $|H \cap D|$ such simple closed curves.
\end{proof}

In short, if $n$ is sufficiently large then $(H + nK) \cap D$ cuts $D$
into pieces and most of these pieces are the parallel annuli, $A_i$.

\subsection{Finding a ``cover'' of K}
%%% Really should be $K$ but hyperref complains.

Recall that $K \setminus \bdy D = K' \disjoint K''$.  Let
$\{\alpha_j'\} = H \cap K'$.  Similarly let $\{\alpha_j''\} = H \cap
K''$.  Due to the minimality assumptions (see the beginning of
Section~\ref{Sec:Assumptions}) every loop of $H \cap K$ is essential
in $K$ and every arc $\alpha_j' \subset K'$ and $\alpha_j'' \subset
K''$ is also essential

Choose a collection of oriented arcs $\{\beta_j'\}$ with the following
properties:

\begin{itemize}
\item
Every arc $\beta_j'$ is simple and is embedded in $K'$.
\item
Both endpoints of $\beta_j'$ are at the point $x$.
\item
The interiors of the $\beta_j'$ are disjoint.
\item
The union of the $\beta_j'$, together with $\bdy D$, forms a
one-vertex triangulation of $K'$.
\item
The chosen arcs $\{\beta_j'\}$ minimizes the quantity 
$|(\bigcup_j \alpha_j') \cap (\bigcup_j \beta_j')|$.
\end{itemize}

%%% This forces the \alpha_j to be normal.

Similarly choose a collection of arcs $\{\beta_j''\}$ for $K''$.  

Now lift everything to a subsurface of $H + nK$ which is ``almost''
a cyclic cover of $K$:  Let $\lift{K} = (H + nK) \cap \neigh(K)$.  Let
$\pi \from \lift{K} \to K$ be the natural projection map.  So $\pi$ is
the composition of the homeomorphism of $\neigh(K) \homeo K \cross (0,
1)$ with projection onto the first factor, restricted to $\lift{K}
\subset \neigh(K)$.  (It is necessary to slightly tilt the vertical
annuli coming from $H \setminus nK$. This makes $\pi$ a local
homeomorphism.)

Thus $\{x_i\} = \pi^{-1}(x)$.  As discussed above for most values
of $i$ the curve $\gamma_i$ is the component of $\pi^{-1}(\bdy D)$
which contains $x_i$.

\begin{figure}[ht]
\psfrag{alpha1}{$\alpha_{j,1}'$}
\psfrag{alpha2}{$\alpha_{j,2}'$}
\psfrag{alpha3}{$\alpha_{j,3}'$}
\psfrag{H}{$H$}
\psfrag{H+nK}{$H + nK$}
$$\begin{array}{c}
\epsfig{file=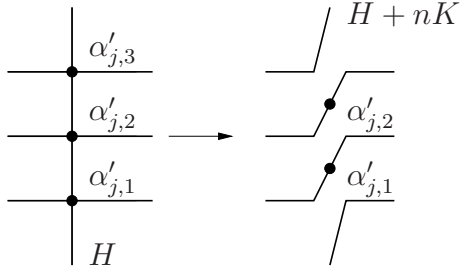, height = 3.5 cm}
\end{array}$$
\caption{The left is before the Haken sum and the right is after.  We
  have tilted the vertical annuli of $H$.}
\label{Fig:Alphas}
\end{figure}

Now lift the set of dual curves $\alpha', \alpha'', \beta', \beta''$:
To be precise, let $\alpha_{j,i}'$ be the component of
$\pi^{-1}(\alpha_j')$ which is contained in the annulus connecting
$K_i$ and $K_{i+1}$.  Define $\alpha_{j,i}''$ similarly.  See
Figure~\ref{Fig:Alphas}.

Let $\beta_{j,i}'$ to be the component of $\pi^{-1}(\beta_j')$ which,
given the orientation of $\beta_j'$, starts at the point $x_i$.
Define $\beta_{j,i}''$ similarly.  Not every $\beta_{j,i}'$ is useful.
However, letting
\begin{equation}
\label{Eqn:AlphaCapBeta}
\alphacapbeta = \max_{k} \left\{\left| \left(\bigcup_j
      \alpha_j'\right) \cap \beta_k' \right|, \left| \left(\bigcup_j
      \alpha_j''\right) \cap \beta_k'' \right| \right\}
\end{equation}
we have:

\begin{claim}
\label{Clm:Betas}
For all $i$ with $\alphacapbeta < i < n - \alphacapbeta$ and for all
$j$ the map $\pi|\beta_{j, i}'$ is onto its image.  The same holds
for $\pi|\beta_{j, i}''$.
\end{claim}

\begin{proof}
Every time $\beta_{j,i}'$ crosses one of the $\alpha_{j,i}'$'s it goes
up (or down) exactly one level.  Thus any $\beta_{j,i}'$, with $i$ as
in the hypothesis, has both endpoints on some lift of $x$ and the
claim holds.
\end{proof}

\begin{define}
\label{Def:BetaShifts}
Suppose that $\alphacapbeta < i < n - \alphacapbeta$.  Suppose that
the final point of $\beta_{j.i}'$ is $x_k$.  By definition of
$\beta_{j,i}'$ the starting point is $x_i$.  Define the {\em shift} of
$\beta_{j,i}'$ to be $\sigma(\beta_{j.i}') = k - i$.
\end{define}

An important observation is:

\begin{claim}
\label{Clm:Shift}
The shift $\sigma(\beta_{j,i}')$ does not depend on the value of $i$.
\qed
\end{claim} 

Henceforth we will use $\sigma(\beta_j')$ to denote the shift of
$\beta_{j,i}'$, for any $i$.  The same notation will be used for arcs
of $K''$.

\subsection{Finding essential curves and annuli}

Now to gain some control over the parallel curves $\gamma_i
\subset D$:  Let 
\begin{equation}
\label{Eqn:LcmNumber}
\lcmnumber = \min_{j,k} \{\lcm(|\sigma(\beta_j')|, |\sigma(\beta_k''))|\}.  
\end{equation}
Here we adopt two conventions: First, the least common multiple of any
number with zero is $\infty$ and second, the minimum of the set
$\{\infty\}$ is zero.  As a consequence, if all shifts are zero in
either $K'$ or $K''$ then $\lcmnumber = 0$.  Finally set
\begin{equation}
\label{Eqn:Enormous}
\enormous = \max \left\{ \HcapBdyD, \alphacapbeta, \lcmnumber \right\}.
\end{equation}

\begin{claim}
\label{Clm:Gammas}
For all $i$ with $\enormous < i < n - \enormous$ the curve $\gamma_i$ is
essential in $H + nK$.
\end{claim}

\begin{proof}
Consider some curve $\gamma_i$ with $i$ in the indicated range.  

First suppose that all shifts on one side, say $K'$, are zero.  That
is, for all $j$ we have $\sigma(\beta_j') = 0$.  Then it follows that
$(H + nK) \setminus \gamma_i$ has two components one of which is
homeomorphic to $K'$ (and in fact is isotopic, relative to $x_i$, to
$K_i'$).  Recall that $n > 2 \cdot \enormous$, and $\euler(K) < 0$,
and Euler characteristic is additive under Haken sum.  Thus
$\euler(K') + 1 > \euler(H) + n\euler(K) = \euler(H + nK)$.  Now, if
$\gamma_i$ is inessential then the other component of $(H + nK)
\setminus \gamma_i$ is a disk.  So the surface $H + nK$ would be
obtained by gluing a copy of $K'$ to a disk along their common
boundary.  It would follow that $\euler(H + nK) = \euler(K') + 1$, a
contradiction.  So if all shifts on one side are zero then $\gamma_i$
is essential.

Now suppose that there are nonzero shifts on both sides.  Reversing
the orientation of some $\beta_j'$ or $\beta_k''$ we may assume that
the shifts $\sigma(\beta_j') = r$ and $\sigma(\beta_k'') = s$ are both
positive and $t = \lcm(r,s) = \lcmnumber$.  Let $\delta'$ be the union
of paths $\bigcup_{l = 0}^{\frac{t}{r} - 1} \beta_{j,i + rl}'$.  Let
$\delta''$ be the union of paths $\bigcup_{l = 0}^{\frac{t}{s} - 1}
\beta_{k,i + sl}''$.  Then $\delta = \delta' \cup \delta''$ is a
simple closed curve in $\lift{K}$ which meets $\gamma_i$ exactly once
at the point $x_i$.  So $\gamma_i$ is essential.
\end{proof}

Similar ideas will give some control over the annuli $A_i \subset D$.
Recall that $\bdy A_i = \gamma_i \cup \gamma_{i+1}$.

\begin{claim}
\label{Clm:Annuli}
For all $i$ with $\enormous < i < n - \enormous - 1$ the annulus $A_i$ is not
boundary parallel into $H + nK$.
\end{claim}

\begin{proof}
Suppose that $A_i$ is boundary parallel into $H + nK$.  Let $B \subset
H + nK$ be the annulus with which $A_i$ cobounds a solid torus.
%%% interior disjoint from H + nK
As the other case is identical, suppose that $B$ is adjacent to the
curves $\gamma_i$ and $\gamma_{i+1}$ from the $K'$-side.  Now, by
Claim~\ref{Clm:Betas}, all of the $\beta_{j,i}'$ and all of the
$\beta_{j,i+1}'$ lie in $B \subset \lift{K} = (H + nK) \cap
\neigh(K)$.  Since $\bdy B = \bdy A_i = \gamma_i \cup \gamma_{i+1}$ all
shifts of the $\beta_{j,i}'$ are zero or one.  Likewise all shifts of
the $\beta_{j,i+1}'$ are zero or minus one.  Thus there is an $i$ with
$\sigma(\beta_{j,i}') = 1$ but $\sigma(\beta_{j,i+1}')$ equals zero
or minus one.  This contradicts Claim~\ref{Clm:Shift}.
\end{proof}

\subsection{Finishing the proof of the theorem}
Recall that all of the curves $\gamma_i$ bound embedded disks in the
manifold because they bound disks in $D$.  Thus by Scharlemann's ``No
Nesting'' Lemma~\ref{Lem:NoNesting}, all of the $\gamma_i$'s bound
disks in one of the two handlebodies bounded by $H + nK$, $V_n$ or
$W_n$.  From strong irreducibility of $H + nK$ and
Claim~\ref{Clm:Gammas} it follows that all the $\gamma_i$'s bound
essential disks on the same side.  As the other case is identical,
suppose that $\gamma_i$ bounds $D_i \subset V_n$ for all $i$.

Now either $A_i$ or $A_{i+1}$ lies in the opposite handlebody $W_n$.
As the two possibilities are symmetric, suppose $A_i \subset W_n$.
There are two final cases.  If $A_i$ is compressible in $W_n$ then
compress to obtain two disks, say $E_i, E_{i+1} \subset W_n$.  Here
$\bdy E_i = \gamma_i = \bdy D_i$.  It follows that $H + nK$ is
reducible, a contradiction.

Suppose instead that $A_i$ is incompressible.  Since $A_i$ is not
boundary parallel (Claim~\ref{Clm:Annuli}) there is a boundary
compression of $A_i$ yielding an essential disk $E_i$ with $\bdy E_i$
disjoint from $\bdy A_i = \gamma_i \cup \gamma_{i+1}$.  So $H + nK$ is
weakly reducible, another contradiction.  This final contradiction
completes the proof of Theorem~\ref{Thm:HPlusNK}. \qed

\section{Adding copies of a torus}
\label{Sec:HPlusNT}

For the remaining part of Theorem~\ref{Thm:HPlusNKGeneral} the surface
added is a torus, $T$.  Hence we deal with sequences of strongly
irreducible Heegaard splittings of the form $H + nT$.

\begin{theorem}
\label{Thm:HPlusNT}
Suppose $M$ is a closed, orientable three-manifold and $H$ and $T$ are
closed orientable transverse surfaces in $M$, with $T$ a two-torus.
Suppose that a Haken sum $H + T$ is given so that the surface $H + nT$
is a strongly irreducible Heegaard splitting for arbitrarily large
values of $n$.  Assume also that no pair of these splittings are
isotopic in $M$.  Then the surface $T$ is incompressible.
\end{theorem}

Assume that $T$ is compressible to obtain a contradiction.  As $M$ is
irreducible there are two cases: Either $T$ bounds a solid torus or
$T$ bounds a {\em cube with a knotted hole}.  Denote the submanifold
which $T$ bounds by $X \subset M$.

Before considering these cases in detail, apply
Lemma~\ref{Lem:RemoveTrivial} so that $H \cap T$ consists of curves
essential on $T$.  These all have the same slope.  Further, assign a
parity to the curves of $H \cap T$ as follows: Choose any oriented
curve $\alpha$ in $T$ which meets each of the components of $H \cap T$
exactly once.  Then, travelling along $\alpha$ in the chosen direction
we cross the curves of $H \cap T$ and, according to the Haken sum, $H
+ nT$ either descends into the submanifold $X$ or ascends out of $X$.
Assign the former a negative parity and the latter a positive.  As the
other case is similar, we assume that there are more curves of $H \cap
T$ of positive parity than negative.  (There cannot be equal numbers
of both as then, for large values of $n$, the surface $H + nT$ fails
to be connected.)  We now have:

\pagebreak

\begin{lemma}
\label{Lem:RemoveCancelling}
Suppose the sequence $H + nT$ is essential in $T$.  Let $m$ be the
number of positive curves of $H \cap T$ minus the number of negative.
Let $m' = (|H \cap T| - m)/2$.  Then there is an isotopy of $H' = H +
m'T$ so that
\begin{itemize}
\item
all curves of $H' \cap T$ are essential in $T$,
\item 
all curves of $H' \cap T$ are positive, and 
\item
for all $n > m'$ the surface $H + nT$ is isotopic to $H' + (n-m')T$.
\end{itemize} \qed
\end{lemma}

As the proof of Lemma~\ref{Lem:RemoveCancelling} is essentially
identical to that of Lemma~\ref{Lem:RemoveTrivial} we omit it.
An essential sequence $H + nT$ {\em reduced} if all of the
curves of $H \cap T$ have the same parity.

\subsection{Bounding a solid torus}

Suppose now that $T$ bounds a solid torus $X$.  We have:

\begin{claim}
\label{Clm:IsotopingBack}
If $H + nT$ is reduced and $m = |H \cap T|$ then, for any positive $n$,
the surface $H + nT$ is isotopic in $M$ to $H + (n + m)T$.
\end{claim}

\begin{proof}
Choose a homeomorphism $X \homeo \DD^2 \cross S^1$, where
$\closure{\neigh(T) \cap X} \homeo A \cross S^1$ with $A \homeo \{z
\in \CC \st 1/2 \leq |z| \leq 1\}$.  Set $D_0 = \closure{\DD^2
\setminus A}$.  

If the slope of $H \cap T$ is meridional (isotopic to $\bdy \DD^2
\cross \{\pt\}$) then the desired isotopy is $\varphi \from M \cross I
\to M$ with $\varphi_t|(M \setminus X) = \identity$, $\varphi_t(z,
\theta) = (z, \theta \pm 2t\pi)$ for all $z \in D_0$, and
$\varphi_t(z, \theta) = (z, \theta \pm 2t\pi \cdot (2 - 2|z|))$ for
all $z \in A$.  Here the sign $\pm$ is determined by the parity of the
curves $H \cap T$.  Note also that we only need to do this isotopy
once, not $m$ times.

For any other slope the desired isotopy is $\varphi \from M \cross I
\to M$ with $\varphi_t|(M \setminus X) = \identity$, $\varphi_t(z,
\theta) = (z \cdot \exp(\pm 2t\pi i), \theta)$ for all $z \in D_0$,
and $\varphi_t(z, \theta) = (z \cdot \exp(\pm 2t \pi i (2 - 2|z|)),
\theta)$ for all $z \in A$.  Again the sign $\pm$ is determined by the
parity of the curves $H \cap T$.
\end{proof}

Thus, when $T$ bounds $X$ a solid torus, the sequence $H + nT$
contains only finitely many isotopy classes of Heegaard splittings.
This is a contradiction.

\subsection{Bounding a cube with a knotted hole}

Suppose now that the two-torus $T$ bounds $X$ a {\em cube with a
knotted hole}.  That is, $X \subset M$ is a submanifold contained in a
three-ball $Y \subset M$, and $T = \bdy X$ compresses in $Y$ but not
in $X$.  The unique slope of this compressing disk is called the {\em
meridian}.

We require one more definition: A pair of transverse surfaces $H$ and
$K$ in a three-manifold $M$ are {\em compression-free} if all curves
of $H \cap K$ are essential on both surfaces.

The main theorem of~\cite{KobayashiRieck04} is: 

\begin{theorem}
\label{Thm:LocalDetectionInCubeWithKnottedHole}
Suppose $H \subset M$ is strongly irreducible and the two-torus $T$
bounds $X \subset M$, a cube with a knotted hole.  Suppose also that
$H$ and $T$ are compression-free with non-trivial intersection.  Then:
\begin{itemize}
\item
the components of $H \cap X$ are all annuli and
\item
there is at least one component of $H \setminus T$ which is an
meridional annulus, boundary parallel into $T$.
\end{itemize}
\end{theorem}

So, choose $H$ and $T$ as provided by the hypotheses of
Theorem~\ref{Thm:HPlusNT}.  Suppose also, as provided by
Lemmas~\ref{Lem:RemoveTrivial} and~\ref{Lem:RemoveCancelling}, that $H
+ nT$ is reduced -- all curves of $H \cap T$ are essential and of the
same parity.

\begin{claim}
\label{Clm:MeridionalIntersection}
All curves of $H \cap T$ are meridional on $T$.
\end{claim}

\begin{proof}
If $H$ and $T$ are compression-free then apply
Theorem~\ref{Thm:LocalDetectionInCubeWithKnottedHole} and we are done.  If
not then there is a curve of intersection which bounds an innermost
disk in $H$ and which is essential on $T$.  As $T$ is not compressible
into $X$ we are done.
\end{proof}

The proof of Theorem~\ref{Thm:HPlusNT}, with $X$ a cube with knotted
hole, now splits into two subcases.  Either $H \cap T$ is
compression-free or not.

\subsubsection{The compression-free case}

Suppose that $H \cap T$ is compression-free and that $H + nT$ is a
reduced sequence.  We again wish to prove that infinitely many of the
$H + nT$ are pairwise isotopic.

Take $nT$ to be $n$ parallel copies of $T$, all inside of $X$.  Note
that $H \cap T = (H + nT) \cap T$ and $H \setminus X = (H + nT)
\setminus X$.  Hence $H + nT$ and $T$ are compression-free.
%%% No disk patches!

We repeatedly isotope $H + nT$ via the following procedure: Apply
Theorem~\ref{Thm:LocalDetectionInCubeWithKnottedHole} to $H + nT$ and
$T$.  Thus there is a meridional annulus $A \subset (H + nT) \setminus
T$ which is boundary parallel into $T$.  Let $B \subset T$ be the
annulus to which $A$ is parallel.  Denote by $Z$ the solid torus which
$A$ and $B$ cobound.

Now, if $A \subset M \setminus X$ then $Z \cap X = B$.  In this case
isotope $A$ and all components of $(H + nT) \cap Z$ into $X$.  Begin
the procedure again applied to this new position of $H + nT$.  

\begin{figure}[ht]
\psfrag{X}{$X$}
\psfrag{Z}{$A$}
\psfrag{T}{$T$}
\psfrag{H+nT}{$H + nT$}
$$\begin{array}{c}
\epsfig{file=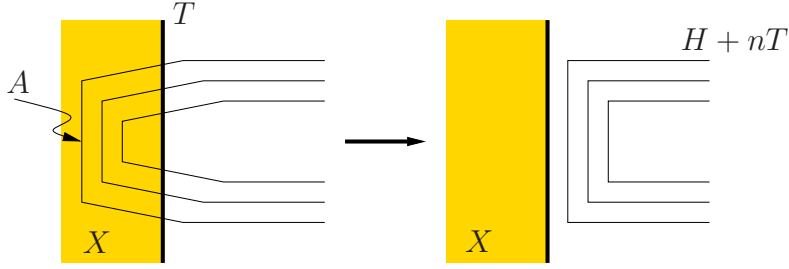, height = 3.5 cm}
\end{array}$$
\caption{Isotoping pieces of $H + nT$ out of $X$.}
\label{Fig:AnnulusIsotopy}
\end{figure}

If $A \subset X$ then $Z \subset X$ as well.  In this case all
components of $(H + nT) \cap Z$ are meridional annuli which are
parallel rel boundary into $T$.  Isotope $A$ and all of the annuli of
$(H + nT) \cap Z$ out of $X$, but keeping them parallel to $T$.  See
Figure~\ref{Fig:AnnulusIsotopy}.

At the end of the procedure, we have isotoped $H + nT$ out of $X$.
The surface $H + nT$ is thus isotopic to a surface which is a union of
components of $H \setminus X$ together with a union of annuli parallel
to sub-annuli of $T$.  There are only finitely many of the latter (as
$H \cap T$ is bounded).  This is a contradiction.

\subsubsection{The meridional compression case}

Suppose now that $H \setminus X$ contains a meridional disk $D \subset
H$ for $T$.  Let $Y$ be the three-ball $X \cup \neigh(D)$.  Note
that all the curves $\{\gamma_j\} = H \cap \bdy Y$ are parallel in
$\bdy Y$.  This is because all of the curves $(H + nT) \cap T$ are
meridional for $T$.  We think of $Y$ as a copy of $\DD^2 \cross I$
-- ``a tall tuna can'' -- with all of the $\gamma_j$ of the form
$\bdy\DD^2 \cross \{\pt\}$.  

\begin{figure}[ht]
\psfrag{D}{$D^i$}
\psfrag{Y}{$Y^i$}
\psfrag{H}{$H^i \setminus Y^i$}
$$\begin{array}{c}
\epsfig{file=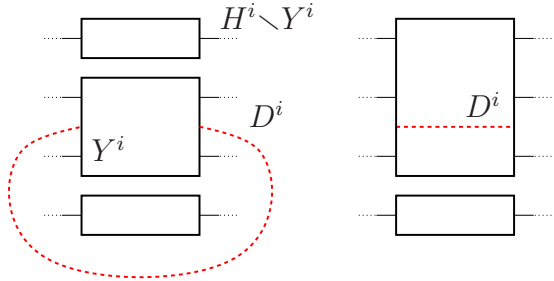, height = 3.5 cm}
\end{array}$$
\caption{The packing step is illustrated on the left while the slicing step is
  on the right.  The disk $D^i$ is depicted by the dotted line.}
\label{Fig:PackingSlicing}
\end{figure}

For each $n$ we carry out an inductive procedure: Fix $n$.  Let $Y^0
= Y$ and let $H^0 = H_n = H + nT$.  At stage $i$ there is a ``stack
of tuna cans'' $Y^i \homeo \DD^2 \cross I_i \subset Y^0$ where $I_i$
is a disjoint union of finitely many closed intervals in $I$.  See
either side of Figure~\ref{Fig:PackingSlicing}.

Each component of $\bdy Y^i$ contains at least one of the curves
$\gamma_j$.  Also, the surface $H^0$ has been isotoped to a surface
$H^i$ so that $H^i \setminus Y^i \subset H^{i-1} \setminus Y^{i-1}
\subset H \setminus Y$.  It follows that $\bdy Y^i \cap H^i$ is a
subset of $\cup \gamma_j$.  Note that all components of $\bdy Y^i
\setminus H^i$ are ``vertical'' annuli or disks.

Suppose some annulus component of $\bdy Y^i \setminus H^i$ is
compressible in $M \setminus (Y^i \cup H^i)$.  So do the ``packing
tuna'' isotopy: There is an disk $D^i$ with interior in $M \setminus
(Y^i \cup H^i)$ and with boundary $\bdy D^i \subset \bdy Y^i$ (see
left side of Figure~\ref{Fig:PackingSlicing}).  Let $Z$ be the
component of $Y^i$ containing $\bdy D^i$.  Then $\bdy D^i$ bounds two
disks in $\bdy Z$, say $E$ and $E'$.  Then either $D^i \cup E$ or $D^i
\cup E'$ bound a three-ball in $M$ which has interior disjoint from
$Z$.  (This is because $M$ is irreducible.)  So there is an isotopy of
$H^i$ which moves some components of $H^i \setminus Z$ into $Z$.  This
reduces the number of curves of intersection $H^i \cap \bdy Y^i$.  Let
$H^{i+1}$ be the new position of $H + nT$.  Let $Y^{i+1}$ be equal to
the union of all components of $Y^i$ which meet $H^{i+1}$.  The
induction hypotheses clearly hold.

Suppose instead some annulus component of $\bdy Y^i \setminus H^i$ is
compressible in $Y^i \setminus H^i$.  Next, perform the ``slice a can in
half'' move: Let $D^i \subset Y^i$ be such a compressing disk with
$\bdy D^i = \DD^2 \cross \{\pt\}$ and $D^i \cap H^i = \emptyset$.  See
right side of Figure~\ref{Fig:PackingSlicing}.  Isotope $D^i \cup H^i$
until $D^i$ is level ($D^i = \DD^2 \cross \{\pt\}$) while maintaining
$D^i \cap H^i = \emptyset$.  This isotopy is supported inside of
$Y^i$.  Let $H^{i+1}$ be the new position of $H + nT$ and let $Y^{i+1}
= Y^i \setminus \neigh(D^i)$.  Again the induction hypotheses clearly
hold.

The procedure terminates after at most $|\{\gamma_j\}| = |(H + nT)
\cap Y|$ steps.  To see this, note that we can never have $|Y^i|$
greater than the original number of curves $\{\gamma_j\}$.  So we
cannot ``slice'' more than that number of times.  Also, the number of
components of $(H + nT) \setminus Y = H \setminus Y$ is bounded and
$H^i \setminus Y^i$ is contained in $H \setminus Y$.  So we cannot
``pack'' more than that number of times.

Let $\last$ be the largest value of $i$ reached in the above
procedure.  After the procedure terminates we have every component of
$\bdy Y^\last \setminus H^\last$ being incompressible in both $M
\setminus (Y^\last \cup H^\last)$ and inside $Y^\last \setminus
H^\last$.  An innermost disk argument shows that every component of
$\bdy Y^\last \setminus H^\last$ is incompressible in $M \setminus
H^\last$.

Let $Z$ be a component of $Y^\last$.  Recall that the curves $\gamma_j
\subset \bdy Z$ are parallel.  Now apply Scharlemann's Local Detection
Theorem~\cite{Scharlemann98} (for three-balls) to $\bdy Z$.  It
follows that $H^\last \cap Z$ is either a disk or is an unknotted
annulus.  

At the end of the procedure the surface $H + nT$ has been isotoped to
a surface which is a union of components of $H \setminus Y$ together
with a union of ``vertical'' annuli and disks of the form $\DD^2
\cross \{\pt\}$.  There are only finitely many of the latter (as $H
\cap \bdy Y$ is bounded).  So for all $n$ the splitting $H+nT$ is
isotopic to one of these finitely many surfaces, a contradiction.
This completes the proof of Theorem~\ref{Thm:HPlusNT}.  \qed

\section{New examples}
\label{Sec:NewExample}

The goal of the next two sections is to give new examples of $H, K,
H+K \subset M$ such that for all integers $n$ the surface $H + nK$ is
a strongly irreducible Heegaard splitting.  

Note that the manifolds of Casson-Gordon have Heegaard genus four and
larger.  Our examples have genus as low as three.  Also, our examples,
unlike those of~\cite{Kobayashi92} and~\cite{LustigMoriah00}, do not
involve twisting around a two-sphere in $S^3$ or require the existence
of an incompressible spanning surface.

In the next two sections we first (\ref{Sec:NewExampleConstruction})
construct our new examples and then (\ref{Sec:NewExampleIsCorrect})
prove that they have the desired properties.

\subsection{Constructing the new examples}
\label{Sec:NewExampleConstruction}

To begin with we sketch the construction, which has obvious
generalizations.  Take $V$ a handlebody of genus three or more.  Take
$\gamma$ to be a ``sufficiently complicated'' curve in $H = \bdy V$.
Double $V$ across $H$ and let $W$ be the other copy of $V$.  Alter the
gluing of $V$ to $W$ by Dehn twisting along $\gamma$ at least five
times.  This gives $M$, a closed orientable manifold.  Now, we will
have a properly embedded surface $K' \subset V$ with $K' \cap \gamma =
\emptyset$.  Thus $K'$ doubles to give a surface $K$ in $M$.  Adding
copies of $K$ to $H$ will give the desired sequence of Heegaard
splittings.

Before giving the details recall:

\begin{define}
\label{Def:DiskBusting}
Let $V$ be a handlebody.  A simple closed curve $\gamma \subset \bdy
V$ is {\em disk-busting} if $\bdy V \setminus \gamma$ is
incompressible in $V$.
\end{define}

\begin{figure}[ht]
\psfrag{g'}{$\gamma'$}
$$\begin{array}{c}
\epsfig{file=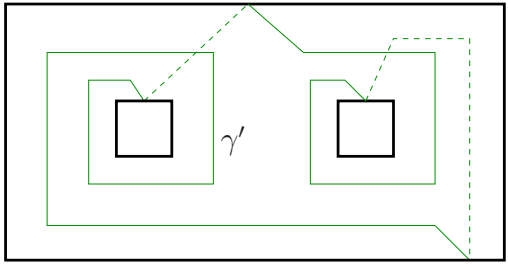, height = 3.5 cm}
\end{array}$$
\caption{The curve $\gamma'$ is disk-busting in $V'$}
\label{Fig:DiskBusting}
\end{figure}

For the remainder of this section take $V'$ a handlebody of genus two.
(Larger genus is also possible.)  Let $\gamma' \subset V'$ be a
non-separating disk-busting curve.  Set $K' = \bdy V' \setminus
\neigh(\gamma')$.  For an example of this see
Figure~\ref{Fig:DiskBusting}.

%%% Fix overcrossings in general.

Take $U$, a solid torus, and fix a subdisk of the boundary $E \subset
\bdy U$.  Let $V'' = (K' \cross I) \cup U$ where $K' \cross I$ is
glued to $U$ via some homeomorphism between a subdisk of $K' \cross
\{1\}$ and the disk $E$.  Thus $E$ and any meridional disk of $U$
(which is disjoint from $E$) are essential disks in $V''$.  Let
$\bdy_+ V'' = \closure{ \left( (K' \cross \{1\}) \cup \bdy U \right)
\setminus E }$.  Let $\bdy_- V'' = K' \cross \{0\}$.

\begin{figure}[ht]
\psfrag{g}{$\gamma$}
$$\begin{array}{c}
\epsfig{file=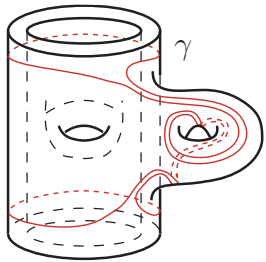, height = 3.5 cm}
\end{array}$$
\caption{The curve $\gamma \subset \bdy_+ V''$ is disk-busting for $V''$.}
\label{Fig:DiskBustingTwo}
\end{figure}

Now choose $\gamma \subset \bdy_+ V''$ a disk-busting curve for $V''$.
See Figure~\ref{Fig:DiskBustingTwo}, for example.

\begin{figure}[ht]
\psfrag{g}{$\gamma$}
$$\begin{array}{c}
\epsfig{file=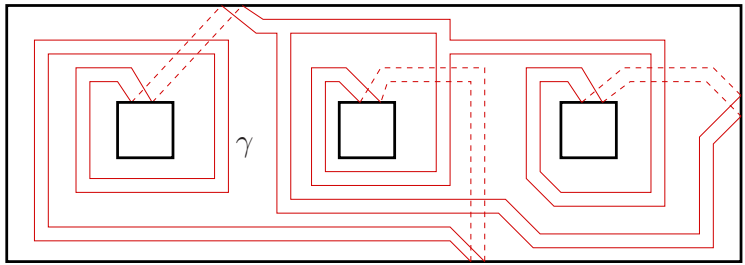, height = 3.5 cm}
\end{array}$$
\caption{To obtain $M$, double the handlebody shown and Dehn twist at
  least five times along $\gamma$.}
\label{Fig:GammaInV}
\end{figure}

Form a genus three handlebody $V$ by gluing $V'$ to $V''$ via the
natural map between $K' \subset \bdy V'$ and $\bdy_- V'' \subset \bdy
V''$.  It is easy to check that $\gamma$ is disk-busting in $V$.  As
this fact is not needed in what follows we omit the proof.  However,
see Figure~\ref{Fig:GammaInV} for a picture.

Now, form a manifold $D(V)$ by {\em doubling $V$} -- that is, let
$W$ be an identical copy of $V$ and glue these two handlebodies by the
identity map between their boundaries.  Finally, obtain a closed
three-manifold $M$ by altering the gluing between $V$ and $W$ by Dehn
twisting at least five times along $\gamma$.  Again, we do not need
the fact that $H$ is a strongly irreducible Heegaard splitting, nor
the consequence that $M$ is irreducible.

Let $K = D(K') \subset M$ be the double of $K'$.  As $K'$ is
connected, so is $K$.  The surface $K$ is also incompressible in $M$,
but as this fact is not required in the sequel, we omit any direct
proof.

\begin{figure}[ht]
\psfrag{A}{$A$}
\psfrag{K}{$K$}
\psfrag{H-A}{$H \setminus A$}
\psfrag{H+K}{$H+K$}
\psfrag{sum}{sum}
$$\begin{array}{c}
\epsfig{file=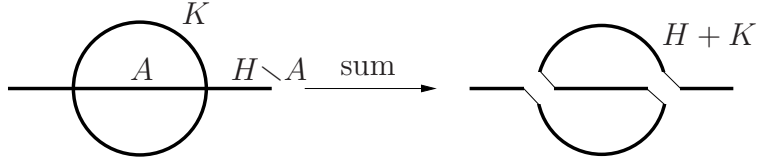, height = 2 cm}
\end{array}$$
\caption{Picture showing (schematically) the relative positions of
  $H$, $K$, and $H + K$.}
\label{Fig:PositionOfHAndK}
\end{figure}

Next, choose the Haken sum of $H$ and $K$: Label the two curves of $K
\cap H$ by $\alpha$ and $\beta$.  Recall that $\gamma'$ was chosen to
be disk-busting and non-separating in $\bdy V'$.  Note that $\alpha$
and $\beta$ cobound an annulus $A = \neigh(\gamma') = \bdy V'
\setminus K' \subset H$ and that $\alpha$ and $\beta$ cut $K$ into two
halves $K' \subset V$ and $K'' \subset W$.  Also, $\alpha$ and $\beta$
cut $H$ into two connected pieces, $A$ and $\closure{H \setminus A}
\homeo \bdy_+ V''$.  Note that $K$ and $H$ are both separating
surfaces in $M$.  For a schematic picture, see the left side of
Figure~\ref{Fig:PositionOfHAndK}.

So choose the Haken sum of $H$ and $K$ as indicated by the right side
of Figure~\ref{Fig:PositionOfHAndK}.  To be precise, let $\mathcal{H}
\from M \cross I \to M$ be an ambient isotopy of $M$ which is fixed
pointwise outside of $\neigh(A)$, moves $\alpha$ across $A$ to
$\beta$, sends the solid torus $\neigh(\alpha)$ to $\neigh(\beta)$,
takes $K \cap \neigh(\alpha)$ to $K \cap \neigh(\beta)$, and takes $H
\cap \neigh(\alpha)$ to $H \cap \neigh(\beta)$.  Now choose any Haken
sum of $H$ and $K$ along $\alpha$ and use $\mathcal{H}$ to transfer
this choice to $\beta$.  Again, see Figure~\ref{Fig:PositionOfHAndK}.
This defines the Haken sum $H + K$ and thus defines $H + nK$.

\subsection{Demonstrating the desired properties}
\label{Sec:NewExampleIsCorrect}

We now can state:

\begin{theorem}
\label{Thm:NewExample}
Given $V$ and $\gamma$ as above, the surface $H + nK$ is a strongly
irreducible Heegaard splitting of $M$, for any even $n > 0$.
\end{theorem}

\begin{remark}
\label{Rem:NewExample}
In fact $H + nK$ is a strongly irreducible Heegaard splitting for any
integer $n$.  We restrict to $n$ positive and even only for notation
convenience.
\end{remark}

\begin{figure}[ht]
\psfrag{g}{$\gamma$}
$$\begin{array}{c}
\epsfig{file=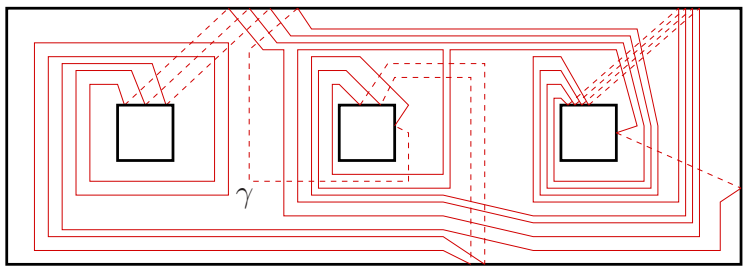, height = 3.5 cm}
\end{array}$$
\caption{Doubling the handlebody and twisting along the curve shown
  gives a hyperbolic manifold satisfying the hypotheses of
  Theorem~\ref{Thm:NewExample}}
\label{Fig:Hyperbolic}
\end{figure}

\begin{remark}
\label{Rem:HyperbolicNewExample}
The curve $\gamma$ in Figure~\ref{Fig:GammaInV} does not give a
hyperbolic manifold because the resulting $M$ contains a pair of Klein
bottles.  See Figure~\ref{Fig:Hyperbolic} for a more complicated curve
$\gamma$.  This curve does yield a hyperbolic manifold with the
desired sequence of Heegaard splittings.
\end{remark}

The proof of Theorem~\ref{Thm:NewExample} divides naturally into two
pieces.

\begin{claim}
\label{Clm:NewExampleIsHeegaardSplitting}
For positive, even $n$ the surface $H + nK$ is a Heegaard splitting.
\end{claim}

\begin{proof}
Recall that $M \setminus \neigh(H \cup K)$ is homeomorphic to the
disjoint union of $V'$, $V''$, $W'$, and $W''$.  Also, the curves $K
\cap H$ are denoted by $\alpha$ and $\beta$.

Let $nK$ be $n$ evenly spaced parallel copies of $K$ in $\neigh(K)$.
That $H + nK$ is connected follows from our choice of Haken sum along
$\alpha$ and $\beta$.  $H + nK$ is separating because $H$ and $K$ are
separating.  See Figure~\ref{Fig:Twisting}.

\begin{figure}[ht]
\psfrag{nK}{$nK$}
\psfrag{H+nK}{$H+nK$}
\psfrag{Vn}{$V_n$}
\psfrag{sum}{sum}
$$\begin{array}{c}
\epsfig{file=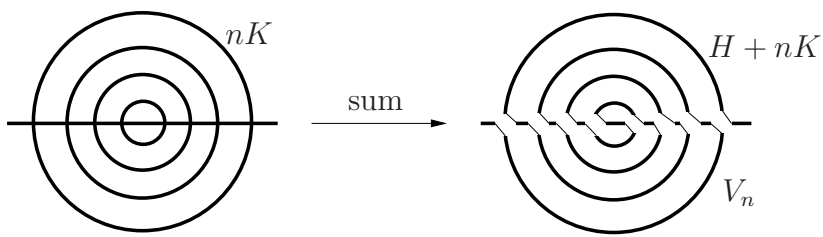, height = 3.0 cm}
\end{array}$$
\caption{Adding $nK$ to $H$ yields a connected, separating surface.}
\label{Fig:Twisting}
\end{figure}

Label the closures of the two components of $M \setminus (H + nK)$ by
$V_n$ and $W_n$ where $V_n$ contains $V \setminus \neigh(K)$ and $W_n$
contains $W \setminus \neigh(K)$. (This is where ``$n$ positive and
even'' is used.  Again, see the right half of
Figure~\ref{Fig:Twisting} for a picture with $n = 4$.)

Consider now the collection of closed annuli $\closure{H \cap
  \interior(V_n)}$.  Cutting $V_n$ along all of these gives several
components: The first, $V'_n$, contains $V' \setminus \neigh(K)$ while
the second, $V''_n$, contains $V'' \setminus \neigh(K)$ and the rest
are isotopic to $\neigh(K')$ or $\neigh(K'')$.  Let $V^P_n$ be the
submanifold of $V_n$ obtained by taking the union of all the latter
(\ie, not $V'_n$ or $V''_n$).  Here the ``$P$'' in the superscript
stands for ``product''.

Let $A_n \cup B_n$ be the two annuli in $\closure{H \cap
  \interior(V_n)}$ which are also in $\bdy V^P_n$.  Here we assign
labels so that $A_n$ meets the component of $H \cap \neigh(K)$ which
contains $\alpha$.  Thus, as $n$ is even, $B_n$ meets the component
of $H \cap \neigh(K)$ which contains $\beta$.  We have realized $V_n$
as the union of three pieces $V'_n$, $V''_n$, and $V^P_n$, glued to
each other along the annuli $A_n$ and $B_n$.

Recall now that $V'_n \homeo V'$, $V''_n \homeo V''$ and thus both
are handlebodies.  Also, the annulus $B_n$ is {\em primitive} in
$V''_n$:  There is a disk in $V''_n$ meeting $B_n$ in a
single co-core arc.  See Figure~\ref{Fig:DiskBustingTwo} and notice
that $B_n$ is parallel to $\beta \cross I \subset \bdy K' \cross I
\subset V''$. 

Since $V^P_n$ and $V''_n$ are handlebodies it follows that $V^P_n
\cup_{B_n} V''_n$ is also a handlebody.  Also, as $V^P_n$ is a
product, the annulus $A_n$ is primitive on $V^P_n \cup_{B_n} V''_n$.
So, since $V'_n$ is a handlebody, we finally have $V_n = V'_n
\cup_{A_n} V^P_n \cup_{B_n} V''_n$ is a handlebody and applying
similar arguments to $W_n$ the surface $H + nK$ is a Heegaard
splitting of $M$.
\end{proof}

\begin{claim}
\label{Clm:NewExampleIsStronglyIrreducible}
For positive, even $n$ the surface $H + nK$ is strongly irreducible.
\end{claim}

\begin{proof}
Recall that $\gamma$ was a curve in $\bdy_+ V''$ and thus also a curve
in $H+nK$.  Recall that $M$ was obtained by doubling $V$ and then
twisting at least five times along $\gamma$.  

We will show that $\gamma$ is disk-busting for $V_n$ and thus for
$W_n$.  The proof of the claim will then conclude with a theorem of
Casson~\cite{MoriahSchultens98} proving that $H + nK$ is strongly
irreducible.

Choose $D$, any essential disk in $V_n$.  Choose a hyperbolic metric
on $H + nK$.  Tighten $\bdy D, \bdy A_n, \bdy B_n, \gamma$ to be
geodesics.  Perform a further isotopy of $D$ relative to $\bdy D$ to
minimize the intersection of $D$ with $A_n \cup B_n$.

Now note that $A_n$ and $B_n$ are incompressible in $V_n$.  If not
then some boundary component of $A_n$ bounds a disk in $V'_n$ or some
boundary component of $B_n$ bounds a disk in $V''_n$.  (None of these
curves bound disks in $V^P_n$ because neither $K'$ nor $K''$ is a
two-sphere.)  The first is impossible because $\bdy A_n$ is parallel
to $\gamma' \subset V'_n$ which is disk-busting.  The second is
impossible because $\bdy_- V''$ is $\pi_1$-injective into $V''_n$.

So no component of $D \cap (A_n \cup B_n)$ is a simple closed curve.
Let $D'$ be an outermost disk of $D \setminus (A_n \cup B_n)$: That
is, $D'$ is the closure of a disk component of $D \setminus (A_n \cup
B_n)$ and $D'$ meets $A_n \cup B_n$ in at most one arc.  It follows
that $D'$ is an essential disk in $V'_n$, $V^P_n$, or $V''_n$.  (If
not we could decrease $|\bdy D \cap (A_n \cup B_n)|$, an
impossibility.)

There are three cases: $D'$ lies in $V'_n$, $V^P_n$, or $V''_n$.

Suppose first that $D' \subset V'_n$.  If $D' = D$ is disjoint from
$A_n$ then, as $A_n$ is parallel to $\gamma'$ in $\bdy V'_n$, we may
isotope $D$ to be disjoint from $\gamma'$.  This contradicts our
choice of $\gamma'$ being disk-busting in $V'_n$.  If $D' \subset D$
is a strict inclusion then $D' \cap A_n$ is a single arc.  Then $D'$
may be isotoped either to lie disjoint from $\gamma'$ ($D' \cap A_n$
is inessential in $A_n$) or to meet $\gamma'$ in a single point ($D'
\cap A_n$ is essential in $A_n$).  Again, this is because $\gamma'$
and $A_n$ are parallel on the boundary on $V'_n$.  The former
contradicts $\gamma'$ being disk-busting.  For the latter take two
parallel copies of $D'$ in $V'_n$ and band these together along
$\gamma' \setminus \neigh(D')$ to obtain an essential disk disjoint
from $\gamma'$.  This is again a contradiction.

The next possibility is that $D'$ lies in $V^P_n$.  However, this
cannot happen as $V^P_n$ is the trivial $I$-bundle over a surface.

We conclude that $D'$ is an essential disk in $V''_n$.  It follows
that $D'$ intersects $\gamma$ because $\gamma$ was chosen to be
disk-busting for $V'' \homeo V''_n$.  Thus $D$ has non-trivial
geometric intersection with $\gamma$.  As our choice of $D$ was
arbitrary we conclude that $\gamma \subset H + nK$ is disk-busting for
both $V_n$ and $W_n$.

Note that $D(V)$, the double of $V$, is reducible.  To obtain $M$ from
$D(V)$ we cut open along a neighborhood of $\gamma$ in $\bdy_+ V''$
and Dehn twisted at least five times.  It follows that $H + nK$ gives
a Heegaard splitting of $D(V)$ and all of these are reducible in
$D(V)$.  (To see this recall that the disk $E$ cut the solid torus $U$
from $V''$.  Thus the double $D(E)$ is a reducing sphere for all of
the $H + nK$ in $D(V)$.)  Thus we are in a position to apply the
following theorem of Casson (see the appendix
of~\cite{MoriahSchultens98}):

\begin{theorem}
\label{Thm:Casson}
Suppose $\gamma \subset H \subset N$ is a curve on a reducible
Heegaard splitting surface of a closed orientable manifold $N$, and
that $H \setminus \gamma$ is incompressible in $N$.  Cutting $N$ open
along a neighborhood of $\gamma$ in $H$ and Dehn-twisting at least
five times gives a strongly irreducible splitting $H'$ of the new
manifold $N'$.
\end{theorem}

It follows that for all positive, even $n$ the splittings $H + nK$ are
strongly irreducible.  We are done.
\end{proof}

Claim~\ref{Clm:NewExampleIsHeegaardSplitting} and
Claim~\ref{Clm:NewExampleIsStronglyIrreducible} together prove
Theorem~\ref{Thm:NewExample}. \qed

\begin{remark}
There is a well-known relationship, due to
Rubinstein~\cite{Rubinstein97} and Stocking~\cite{Stocking00}, between
strongly irreducible splittings and almost normal surfaces.  In
particular, strongly irreducible surfaces should contain a single
place (or ``site'' in Rubinstein's terminology) where the curvature is
highly negative.  This supposedly corresponds to the almost normal
octagon or annulus of the almost normal surface.  In our examples we
find that the subsurface $\bdy_+ V''$ is the distinguished subsurface
of $H + nK$ which presumably contains this special site.
\end{remark}

\pagebreak

\section{Questions}
\label{Sec:Questions}

Recall that Theorem~\ref{Thm:HPlusNK} was originally conjectured by
Sedgwick along with the much stronger:

\vspace{2mm}
\noindent
{\bf Conjecture~\ref{Conj:Sedqwick}.}
{\em 
Let $M$ be a closed, orientable $3$-manifold which contains infinitely
many irreducible Heegaard splittings that are pairwise non-isotopic.
Then $M$ is Haken.
}
\vspace{1mm}

This conjecture may be split, roughly, into two parts.  First we have
the so-called ``Generalized Waldhausen Conjecture'':

\begin{conjecture}
\label{Conj:GeneralizedWaldhausen}
Let $M$ be a closed, orientable $3$-manifold which contains infinitely
many Heegaard splittings, pairwise non-isotopic, all of the same
genus.  Then $M$ is toroidal.
\end{conjecture}

%%% Say what Johannson has done.  

Note that this has been claimed by Jaco and Rubinstein.  However, no
manuscript is available as of the writing of this paper.

The other half of Sedgwick's conjecture deals with splittings of
increasing genus and was the inspiration for our current work:

\begin{conjecture}
\label{Conj:WeakSedqwick}
Let $M$ be a closed, orientable $3$-manifold which contains
irreducible Heegaard splittings of arbitrarily large genus.  Then $M$
is Haken.
\end{conjecture}

We now turn to questions about examples.  In all of the manifolds
listed above, which contain splittings of arbitrarily large genus, the
three-manifold has had Heegaard genus three or higher.  Kobayashi
asks:

\begin{question}
\label{Que:Genus}
Is there an example of a Heegaard genus two manifold which admits
strongly irreducible splittings of arbitrarily large genus?
\end{question}

%%% Kobayashi suggests that the answer to this question is no.

\begin{remark}
\label{Rem:H+nT}
Note that there are examples of toroidal manifolds containing
infinitely many strongly irreducible splittings all of the form $H +
nT$.  Here $H$ is a genus two Heegaard splitting and $T$ is an
incompressible torus; see~\cite{MorimotoSakuma91}.
\end{remark}

Sedgwick, in~\cite{Sedgwick97}, has shown that the
Casson-Gordon examples satisfy the so-called ``Stabilization
Conjecture~\cite{Scharlemann02}''.
%%% recorded as Conjecture 7.4 in the arXiv version -- is this really
%%% the earliest explicit statement of the conjecture?
That is, for any two splittings $H$ and $H'$ obtained from the same
pretzel knot, after stabilizing the higher genus splitting once we may
destabilize to find the lower genus splitting.  Sedgwick's techniques
apply to all of the splittings discussed in
Section~\ref{Sec:OldExamples}.  Kobayashi suggests that the examples
of $H + n K$ given in this paper, after stabilizing twice, should
destabilize about $2n$ times.
%%% have to check that Kobayashi gets credit for this and not somebody
%%% else... I believe that he mentioned somebody else...

\begin{question}
\label{Que:Stabilization}
Does one stabilization suffice?
\end{question}

\bibliographystyle{plain}
\bibliography{bibfile}
\end{document}